\documentclass[12pt]{amsart}
\usepackage{amssymb}
\usepackage{epsf, epsfig}
\usepackage{multirow}
\usepackage{colordvi}
\numberwithin{equation}{section}

\newtheorem{thm}[equation]{Theorem}
\newtheorem{prop}[equation]{Proposition}
\newtheorem{lemma}[equation]{Lemma}
\newtheorem{cor}[equation]{Corollary}

\newtheorem{conj}[equation]{Conjecture}
\newtheorem{example}[equation]{Example}
\newtheorem{remark}[equation]{Remark}
\newtheorem{definition}[equation]{Definition}

\newenvironment{ex}{\begin{example}\rm}{\end{example}}
\newenvironment{rem}{\begin{remark}\rm}{\end{remark}}
\newcounter{FNC}[page]
\def\newfootnote#1{{\addtocounter{FNC}{2}$^\fnsymbol{FNC}$%
     \let\thefootnote\relax\footnotetext{$^\fnsymbol{FNC}$#1}}}

\newcommand{\precdot}{{\prec\!\!\!\cdot\,}}

\newcommand{\Fl}{{{\mathbb F}\ell}}
\newcommand{\Flan}{{{\mathbb F}\ell(\alpha;n)}}

\newcommand{\Gr}{{\rm Gr}}

\newcommand{\C}{{{\mathbb C}}}
\renewcommand{\P}{{{\mathbb P}}}
\newcommand{\R}{{{\mathbb R}}}
\newcommand{\RP}{{\R\P}}

\newcommand{\calM}{{{\mathcal M}}}

\copyrightinfo{}{}
\title[Experimentation in the real Schubert
 calculus]{Experimentation and conjectures in the real Schubert
 calculus for flag manifolds}

\author{Jim Ruffo}
\address{Department of Mathematics\\
         Texas A\&M University\\
         College Station\\
         TX \ 77843\\
         USA}
\email{jruffo@math.tamu.edu}
\urladdr{http://www.math.tamu.edu/\~{}jruffo}
\author{Yuval Sivan}
\address{Department of Mathematics\\
        University of Massachusetts\\
        Lederle Graduate Research Tower\\
        Amherst, MA, 01003\\
        USA}
\email{yuval@student.umass.edu}
\author{Evgenia Soprunova}
\address{Department of Mathematics\\
        University of Massachusetts\\
        Lederle Graduate Research Tower\\
        Amherst, MA, 01003\\
        USA}
\email{esoprun@math.umass.edu}
\urladdr{http://www.math.umass.edu/\~{}esoprun}
\author{Frank Sottile}
\address{Department of Mathematics\\
         Texas A\&M University\\
         College Station\\
         TX \ 77843\\
         USA}
\email{sottile@math.tamu.edu}
\urladdr{http://www.math.tamu.edu/\~{}sottile}

\thanks{Work and computation done at MSRI supported by NSF grant DMS-9810361}
\thanks{Some computations done on computers purchased with NSF SCREMS grant DMS-0079536}
\thanks{Work of Sottile was supported by the Clay Mathematical Institute}
\thanks{This work was supported in part by NSF CAREER grant DMS-0134860}

\begin{document}
\begin{abstract}
 The Shapiro conjecture in the real Schubert calculus, while likely
 true for Grassmannians, fails to hold for flag
 manifolds, but in a very interesting way.
 We give a refinement of the Shapiro conjecture for flag manifolds and present
 massive computational experimentation in support of this refined conjecture.
 We also prove the conjecture in some special cases using discriminants and
 establish relationships between different cases of the conjecture.
\end{abstract}
\maketitle


\section*{Introduction}
 The Shapiro conjecture for Grassmannians~\cite{So00a,SS02} has driven
 progress in enumerative real algebraic geometry~\cite{So03b}, which is the
 study of real solutions to geometric problems. 
 It conjectures that a (zero-dimensional)
 intersection of Schubert subvarieties of a Grassmannian consists
 entirely of real points---if the Schubert subvarieties are given by
 flags osculating a real rational normal curve.
 This particular Schubert intersection problem is quite natural;
 it can be interpreted in terms of real linear series on $\P^1$ with
 prescribed (real) ramification~\cite{EH83,EH87}, real rational
 curves in $\P^n$ with real flexes~\cite{KhS03}, linear systems
 theory~\cite{RS98}, and the Bethe ansatz and Fuchsian equations~\cite{MV04}. 
 The Shapiro conjecture has implications for all these areas.
 Massive computational evidence~\cite{So00a,Ver00} as well as its
 proof by Eremenko and Gabrielov for Grassmannians of codimension 2
 subspaces~\cite{EG02a} give compelling evidence
 for its validity. 
 A local version, that it holds when the Schubert varieties are 
 special (a technical term) and when the points of osculation are
 sufficiently clustered~\cite{So99a}, showed that the special
 Schubert calculus is fully real (such geometric problems can have
 all their solutions real).
 Vakil later used other methods to show that the general Schubert
 calculus on the Grassmannian is fully real.~\cite{Va04}

 The original Shapiro conjecture stated that such an intersection of Schubert
 varieties in a \Red{{\it flag manifold}} would consist entirely of real
 points.
 Unfortunately, this conjecture fails for the first non-trivial
 enumerative problem on a non-Grassmannian flag manifold, but in a very
 interesting way.  
 Failure for flag manifolds was first noted in~\cite[\S 5]{So00a} and
 a more symmetric 
 counterexample was found in~\cite{So00b}, where computer experimentation 
 suggested that the conjecture would hold if the points where the flags
 osculated the rational normal curve satisfied a certain non-crossing condition.
 Further experimentation led to a precise formulation of this refined 
 non-crossing conjecture in~\cite{So03b}.
 That conjecture was only valid for two- and three- step flag manifolds, and the further
 experimentation reported here leads to versions (Conjectures~\ref{C:Main}
 and~\ref{C:New_Main}) for all flag 
 manifolds in which the points of osculation satisfy a monotonicity condition.

 We have systematically investigated the Shapiro conjecture for flag manifolds to gain
 a deeper understanding both of its failure and of our refinement.
 This investigation includes 15.76 gigahertz-years of computer experimentation, 
 theorems relating our conjecture for different enumerative problems,
 and its  proof in some cases using discriminants.
 Recently, our conjecture was proven by Eremenko, Gabrielov,
 Shapiro, and Vainshtein~\cite{EGSV} for manifolds of flags
 consisting of a codimension 2 plane lying on a hyperplane.
 Our experimentation also uncovered some new and interesting phenomena
 in the Schubert calculus of a flag manifold, and it included
 substantial computation in support of the Shapiro conjecture on the
 Grassmannians $\Gr(3,6)$, $\Gr(3,7)$, and $\Gr(4,8)$.

 Our conjecture is concerned with a subclass of Schubert intersection problems.
 Here is one open instance of this conjecture, expressed as a
 system of polynomials in local coordinates for the variety of flags
 $E_2\subset E_3$ in $5$-space, where $\dim E_i=i$.
 Let $t,x_1,\dotsc,x_8$ be indeterminates, and consider the polynomials
\begin{eqnarray*}
   f(t;x) & :=&   \det\left[
   \begin{array}{ccccc}
    1    & 0   & x_1 & x_2 & x_3\\
    0    & 1   & x_4 & x_5 & x_6\\
    t^4  & t^3 & t^2 &  t  &  1\rule{0pt}{14pt}\\
    4t^3 & 3t^2& 2t  &  1  &  0\\
   12t^2 & 3t  &  2  &  0  &  0
   \end{array}\right]\,, \qquad\mbox{and}\vspace{5pt}\\
   g(t;x)& :=&   \det\left[
   \begin{array}{ccccc}
    1 & 0 & x_1 & x_2 & x_3\\
    0 & 1 & x_4 & x_5 & x_6\\
    0 & 0 &  1  & x_7 & x_8\\
    t^4 & t^3 & t^2 & t & 1\rule{0pt}{14pt}\\
    4t^3 & 3t^2 & 2t &1 & 0\\ 
   \end{array}\right]\, .
\end{eqnarray*}

\noindent{\bf Conjecture A.}
{\it
  Let $t_1<t_2<\dotsb<t_8$ be real numbers.
  Then the polynomial system
\begin{eqnarray*}
  \Blue{f(t_1;x)=f(t_2;x)=f(t_3;x)=f(t_4;x)}&=& 0,\quad\mbox{and}\\ 
   \Magenta{g(t_5;x)=g(t_6;x)=g(t_7;x)=g(t_8;x)}&=& 0
\end{eqnarray*}
 has $12$ solutions, and all of them are real.
}\medskip

Evaluating the polynomial $f$ at points $t_i$ preceeding the points at which the 
polynomial $g$ is evaluated is the monotonicity condition.
If we had switched the order of $\Blue{t_4}$ and $\Magenta{t_5}$, 
\[
  \Blue{t_1}\,<\,\Blue{t_2}\,<\,\Blue{t_3}\,<\quad \Magenta{t_5}\ <\ 
  \Blue{t_4}\quad <\,\Magenta{t_6}\,<\,\Magenta{t_7}\,<\,\Magenta{t_8}\,,   
\]
then this would not be monotone.
We computed 400,000 instances of this polynomial
system at different choices of points $t_1<\dotsb<t_8$ (which were monotone),
and each had 12 real solutions. 
In contrast, there were many non-monotone choices of points for which not all
solutions were real, and the minimum number of real solutions that we observe
seems to depend on the combinatorics of the evaluation. 
For example, the system with interlaced points $t_i$ 
\[
  \Blue{f(-8;x)}=\Magenta{g(-4;x)}=
  \Blue{f(-2;x)}=\Magenta{g(-1;x)}=
  \Blue{f(1;x)}=\Magenta{g(2;x)}=
  \Blue{f(4;x)}=\Magenta{g(8;x)}=0
\]
has 12 solutions, \Red{{\em none}} of which are real.
This investigation is summarized in Table~\ref{table:12-flag}.\smallskip

This paper is organized as follows.
In Section 1, we provide background material on flag manifolds, state
the Shapiro Conjecture, and give a geometrically vivid example of its failure.
In Section 2, we give the results of our experimentation, stating our conjectures and
describing some interesting phenomena that we have observed in our data.
The discussion in Section 3 contains theorems about our conjectures, a
generalization of our main conjecture, and proofs of it in some
cases using discriminants. 
Finally, in Section 4 we describe our methods, explain our experimentation and give
a brief guide to our data, all of which and much more is tabulated and available on line
at {\tt www.math.tamu.edu/\~{}sottile/pages/Flags/}.

We thank the Department of Mathematics and Statistics at the
University of Massachusetts at Amherst and the Mathematical Sciences
Research Institute; most of the experimentation underlying our results
was conducted on computers at these institutions.
Funds from the NSF grants DMS-9810361, DMS-0079536, DMS-0070494, and
DMS-0134860 purchased and maintained these computers.
This project began as a vertically integrated research project in the
Summer of 2003.

\section{Background}

\subsection{Basics on flag manifolds}\label{S:basics}
Given positive integers $\alpha:=\{\alpha_1<\dotsb<\alpha_k\}$ with
$\alpha_k<n$,  let $\Flan$ 
be the manifold of flags in $\C^n$ of \Blue{{\it type}} $\alpha$,
\[
   \Flan\ :=\ \{E_\bullet=E_{\alpha_1}\subset E_{\alpha_2}\subset\dotsb\subset
                 E_{\alpha_k}\subset\C^n\mid \dim E_{\alpha_i}=\alpha_i\}\,.
\]
If we set $\alpha_0:=0$, then this algebraic manifold has dimension 
 \[
    \dim(\alpha)\ :=\  \sum_{i=1}^k (n-\alpha_i)(\alpha_i-\alpha_{i-1})\,.
 \]
\Blue{{\it Complete flags}} in $\C^n$ have type $1<2<\dotsb<n{-}1$.

Define $W^{\alpha}\subset S_n$ to be the set of permutations with descents
in $\alpha$,
 \[ 
   W^{\alpha}\ :=\ \{w\in S_n\mid i\not\in\{\alpha_1,\dotsc,\alpha_k\}\Rightarrow  
      w(i)<w(i+1)\}\,.
 \]
We often write permutations as a sequence of their values, omitting commas if
possible. 
Thus $(1,3,2,4,5)=13245$ and $341526$ are permutations in $S_5$ and $S_6$,
respectively. 
Since a permutation $w\in W^\alpha$ is determined by its values before 
its last descent, we need only write its first $\alpha_k$ values.
Thus $132546\in W^{\{2,4\}}$ may be written $1325$.
Lastly, we write $\sigma_i$ for the simple transposition $(i,i{+}1)$.

The positions of flags $E_\bullet$ of type $\alpha$ relative to a fixed
complete flag $F_\bullet$ stratify $\Flan$ into \Blue{{\it Schubert cells}}.
The closure of a Schubert cell is a \Blue{{\it Schubert variety}}.
Permutations $w\in W^{\alpha}$ index Schubert cells $X^\circ_wF_\bullet$ and 
Schubert varieties $X_wF_\bullet$ of $\Flan$.
More precisely, if we set $r_w(i,j):=|\{l\leq i\mid j+w(l)>n\}|$, then 
 \begin{eqnarray}
   X^{\Red{\circ}}_wF_\bullet&=&\{E_\bullet\mid \dim E_{\alpha_i}\cap F_j\Red{=}\nonumber 
       r_w(\alpha_i,j),\ i=1,\dotsc,k,\ j=1,\dotsc,n\},\quad\mbox{and}\\ 
  X_wF_\bullet&=&\{E_\bullet\mid \dim E_{\alpha_i}\cap F_j\Red{\geq} r_w(\alpha_i,j),\
       i=1,\dotsc,k,\ j=1,\dotsc,n\}\,.
 \end{eqnarray}
Flags $E_\bullet$ in $X^\circ_wF_\bullet$ have 
\Blue{{\it position $w$ relative to $F_\bullet$}}.
We will refer to a permutation $w\in W^{\alpha}$ as a 
\Blue{{\em Schubert condition}} on flags of type $\alpha$. 
The Schubert subvariety $X_wF_\bullet$ is irreducible with codimension
$\ell(w):=|\{i<j\mid w(i)>w(j)\}|$ in $\Flan$.

Schubert cells are affine spaces with  
$X^\circ_wF_\bullet\simeq \C^{\dim(\alpha)-\ell(w)}$.
We introduce a convenient set of coordinates for Schubert cells.
Let $\calM_w$ be the set of $\alpha_k\times n$ matrices, some of whose
entries $x_{i,j}$ are fixed: $x_{i,w(i)} = 1$ for $i=1,\dotsc,\alpha_k$
and $x_{i,j}=0$ if 
\begin{center}
$j<w(i)\ $ or $\ w^{-1}(j)<i\ $ or $\ \alpha_l<i<w^{-1}(j)<\alpha_{l+1}$ for some
$l$, 
\end{center}
%
%
and whose remaining $\dim(\alpha)-\ell(w)$ entries give coordinates for $\calM_w$.
For example, if $n=8$, $\alpha=(2,3,6)$, and $w=25\,3\,167$, then 
$\calM_w$ consists of matrices of the form
\[
  \left(\begin{matrix}
      0   &\Red{1}&x_{13}&x_{14}&0&x_{16}&x_{17}&x_{18}\\
      0   &0&  0   &   0  &\Red{1}&x_{26}&x_{27}&x_{28}\\
      0   &0&\Red{1}&x_{34}&0&x_{36}&x_{37}&x_{38}\rule{0pt}{13pt}\\
   \Red{1}&0&  0   &x_{44}&0&  0   &  0   &x_{48}\rule{0pt}{13pt}\\
      0   &0&  0   &   0  &0&\Red{1}&  0   &x_{58}\\
      0   &0&  0   &   0  &0&  0   &\Red{1}&x_{68}
   \end{matrix}\right)\ .
\]

The relation of $\calM_w$ to the Schubert cell $X_w^\circ F_\bullet$ is as follows.
Given a complete flag $F_\bullet$, choose an ordered basis 
$e_1,\dotsc,e_n$ for $\C^n$ corresponding to
the columns of matrices in $\calM_w$ such that 
$F_i$ is the linear span of the last $i$ basis vectors,
$e_{n+1-i},\dotsc,e_{n-1},e_n$. 
Given a matrix $M\in\calM_w$, set $E_{\alpha_i}$ to be the row space of the
first $\alpha_i$ rows of $M$.
Then the flag $E_\bullet$ has type $\alpha$ and lies in the 
Schubert cell $X^\circ_wF_\bullet$,
every flag $E_\bullet\in X^\circ_wF_\bullet$ arises in this way, and the 
association $M\mapsto E_\bullet$ is an algebraic bijection 
between $\calM_w$ and $X^\circ_wF_\bullet$.
This is a flagged version of echelon forms.  See~\cite{Fu97}  
for details and proofs.

Let $\iota$ be the identity permutation.
Then $\calM_{\iota}$ provides local coordinates for $\Flan$ in which the
equations for a Schubert variety are easy to describe.
Note that  
\[
   \dim (E_{\alpha_i}\cap F_j)\ \geq\ r\quad 
    \Longleftrightarrow\quad \mbox{rank}(A)\ \leq\ \alpha_i+j-r\,,
\]
where the matrix $A$ is formed by stacking the first $\alpha_i$ rows of
$\calM_{\iota}$ on top of a $j\times n$ matrix with row span $F_j$.
Algebraically, this rank condition is the vanishing of all minors of
$A$ of size $1{+}\alpha_i{+}j{-}r$.
The polynomials $f$ and $g$ of Example~A from the Introduction arise
in this way.
There $\alpha=\{2,3\}$ and $\calM_\iota$ is the matrix of variables in
the definition of $g$.

Suppose that $\beta$ is a subsequence of $\alpha$.
Then $W^\beta\subset W^{\alpha}$.
Simply forgetting the components of a flag $E_\bullet\in\Flan$ that do not have
dimensions in the sequence $\beta$ gives a flag in $\mathbb{F}\ell(\beta;n)$.
This defines a map
 \[
   \pi\ \colon\ \Flan\ \longrightarrow\ \Fl(\beta;n)
 \]
whose fibres are (products of) flag manifolds.
The inverse image of a Schubert variety $X_wF_\bullet$ of $\mathbb{F}\ell(\beta;n)$
is the Schubert variety  $X_wF_\bullet$ of $\Flan$.

When $\beta=\{b\}$ is a singleton, $\mathbb{F}\ell(\beta;n)$ is the
Grassmannian of $b$-planes in $\C^n$, written $\Gr(b,n)$.
Non-identity permutations in $W^\beta$ have a unique descent at $b$.
A permutation $w$ with a unique descent is \Blue{{\it Grassmannian}}
as the associated Schubert variety $X_wF_\bullet$ (a 
\Blue{{\it Grassmannian Schubert variety}}) is the inverse image of a Schubert
variety in a Grassmannian.

\subsection{The Shapiro Conjecture}\label{S:SC}
A list $(w_1,\dotsc,w_m)$ of permutations in $W^{\alpha}$ is called 
a \Blue{{\it Schubert problem}} if $\ell(w_{1})+\cdots+\ell(w_{m})=\dim(\alpha)$.
Given such a list and complete flags
$F_\bullet^1,\dotsc,F_\bullet^m$, consider the Schubert intersection
 \begin{equation} \label{SchInt}
   X_{w_{1}}F_\bullet^{1} \cap \dotsb  \cap X_{w_{m}}F_\bullet^{m}\,.
 \end{equation}
When the flags $F_\bullet^i$ are in general position, this intersection is
zero-dimensional (in fact transverse by the Kleiman-Bertini
theorem~\cite{Kl74}), and it equals the intersection of the corresponding
Schubert cells.  
In that case, the intersection~\eqref{SchInt} consists of those
flags $E_\bullet$ of type $\alpha$ which have position $w_i$ relative to $F_\bullet^i$, for
each $i=1,\dotsc,m$. 
We call these \Blue{{\em solutions}} to the Schubert intersection
problem~\eqref{SchInt}. 
The number of solutions does 
not depend on the choice of flags (as long as the intersection is transverse)
and we call this number the \Blue{{\it degree}} of the Schubert problem.
This degree may be computed, for example, in the cohomology ring of
the flag manifold $\Flan$.

The Shapiro conjecture concerns the following variant of this
classical enumerative geometric problem:
Which \Red{{\em real}}\/ flags $E_\bullet$ have given position $w_i$ relative to
\Red{{\em real}}\/ flags $F_\bullet^i$, for each $i=1,\dotsc,m$?
In the Shapiro conjecture, the flags $F_\bullet^i$ are not general real
flags, but rather flags osculating a rational normal curve.
Let $\gamma\colon \C \rightarrow \C^n$ be the rational normal curve, 
$\gamma(t):=(1,t,t^2,\ldots,t^{n-1})$ written with respect to 
the ordered basis $e_1,\dotsc,e_n$ for $\C^n$ given above.
The \Blue{{\it osculating flag}} $F_\bullet(t)$ of subspaces to $\gamma$
at the point $\gamma(t)$ is the flag whose $i$-dimensional component is
\[
   F_i(t)\ :=\ \mbox{span}\{\gamma(t),
    \gamma'(t),\ldots,\gamma^{(i-1)}(t)\}\,.
\]
When $t=\infty$, the subspace $F_i(\infty)$ is spanned
by $\{e_{n+1-i},\dotsc,e_n\}$ and $F_\bullet(\infty)$ is the flag used to 
describe the coordinates $\calM_w$. 
If we consider this projectively, $\gamma\colon\P^1\to\P^{n-1}$ is the rational
normal curve and $F_\bullet(t)$ is the flag of subspaces osculating $\gamma$ at
$\gamma(t)$. 

\begin{conj}[B.~Shapiro and M.~Shapiro]\label{C:SC}
  Suppose that $(w_1,\dotsc,w_m)$ is a Schubert problem for flags of type
  $\alpha$.
  If the flags $F_\bullet^1,\dotsc,F_\bullet^m$ osculate the
  rational normal curve at distinct real points, then the 
  intersection~$\ref{SchInt}$ is transverse and consists only of real points. 
\end{conj}

The Shapiro conjecture is concerned with intersections of the form
 \begin{equation}\label{E:shapint}
     X_{w_1}(t_1)\cap X_{w_2}(t_2)\cap 
    \dotsb\cap X_{w_m}(t_m)\,,
 \end{equation}
where we write $X_w(t)$ for $X_wF_\bullet(t)$.
This intersection is an \Blue{{\it instance}} of the Shapiro conjecture for the
Schubert problem $(w_1,\dotsc,w_m)$ at the points
$(t_1,\dotsc,t_m)$. 

Conjecture~\ref{C:SC} dates from around 1995.
Experimental evidence of its validity for Grassmannians was first found
in~\cite{RS98,So97c}.  
This led to a systematic investigation on Grassmannians, both experimentally and 
theoretically in~\cite{So00a}.
There, the conjecture was proven using discriminants for several
(rather small) Schubert problems and 
relationships between the conjecture for different 
Schubert problems were established.
(See also Theorem~2.8 of~\cite{KhS03}.)
For example, if the Shapiro conjecture holds on a Grassmannian for
the Schubert problem consisting only of codimension 1 (simple)
conditions, then it holds for all Schubert problems on that
Grassmannian and on all smaller Grassmannians, if we drop the claim of
transversality. 
More recently, Eremenko and Gabrielov proved the conjecture for
any  Schubert problem on a Grassmannian of codimension
2-planes~\cite{EG02a}. 
Their result is appealingly interpreted as a rational function all of whose 
critical points are real must be real.

The original conjecture was for flag manifolds, but a counterexample was found
and reported in~\cite{So00a}. 
Subsequent experimentation refined this counterexample, and has 
suggested a reformulation of the original conjecture.
We study this refined conjecture and report on
massive computer experimentation (15.76 gigahertz-years) undertaken in 2003 and
2004 at the University of Massachusetts at Amherst, at the MSRI in
2004, and some at Texas A\&M University in 2005.
A byproduct of this experimentation was the discovery of several new and
unusual phenomena, which we will describe through examples.
The first is the smallest possible counterexample to the original
Shapiro conjecture.

\subsection{The Shapiro conjecture is false for flags in 3-space}\label{S:F124} 
We use $\sigma^b$ to indicate that the Schubert condition $\sigma$ is repeated
$b$ times and write $\sigma_i$ for the simple transposition $(i,i{+}1)$.
Then $\bigl(\sigma_2^3,\,\sigma_3^2\bigr)$ is a Schubert problem for flags of
type $\{2,3\}$ in $\C^4$.
For distinct points $s,t,u,v,w\in\RP^1$, consider the Schubert intersection
 \begin{equation}\label{F234Int}
    X_{\sigma_2}(s)\cap X_{\sigma_2}(t)\cap X_{\sigma_2}(u)\  \cap\ 
   X_{\sigma_3}(v)\cap X_{\sigma_3}(w)\,.
 \end{equation}
As flags in projective 3-space, a 
partial flag of type $\{2,3\}$ is a line $\ell$ lying on a plane $H$.
Then $(\ell\subset H)\in X_{\sigma_2}(s)$ if $\ell$ meets the line $\ell(s)$ 
tangent to $\gamma$ at $\gamma(s)$, and 
$(\ell\subset H)\in X_{\sigma_3}(v)$ if $H$ contains the point
$\gamma(v)$ on the rational normal curve $\gamma$.

Suppose that the flag $\ell\subset H$ lies in the
intersection~\eqref{F234Int}.
Then $H$ contains the two points $\gamma(v)$ and $\gamma(w)$, and hence the
secant line $\lambda(v,w)$ that they span.
Since $\ell$ is another line in $H$, $\ell$ meets this secant line
$\lambda(v,w)$.
As $\ell\neq\lambda(v,w)$, it determines $H$ uniquely as the
span of $\ell$ and $\lambda(v,w)$.
In this way, we are reduced to determining the lines $\ell$ which meet the
three tangent lines $\ell(s)$, $\ell(t)$, $\ell(u)$, and the secant line
$\lambda(v,w)$.

\begin{figure}[htb]
 \[
  \begin{picture}(309.6,144)(0,6)
   \put(0,0){\includegraphics[height=158.4pt]{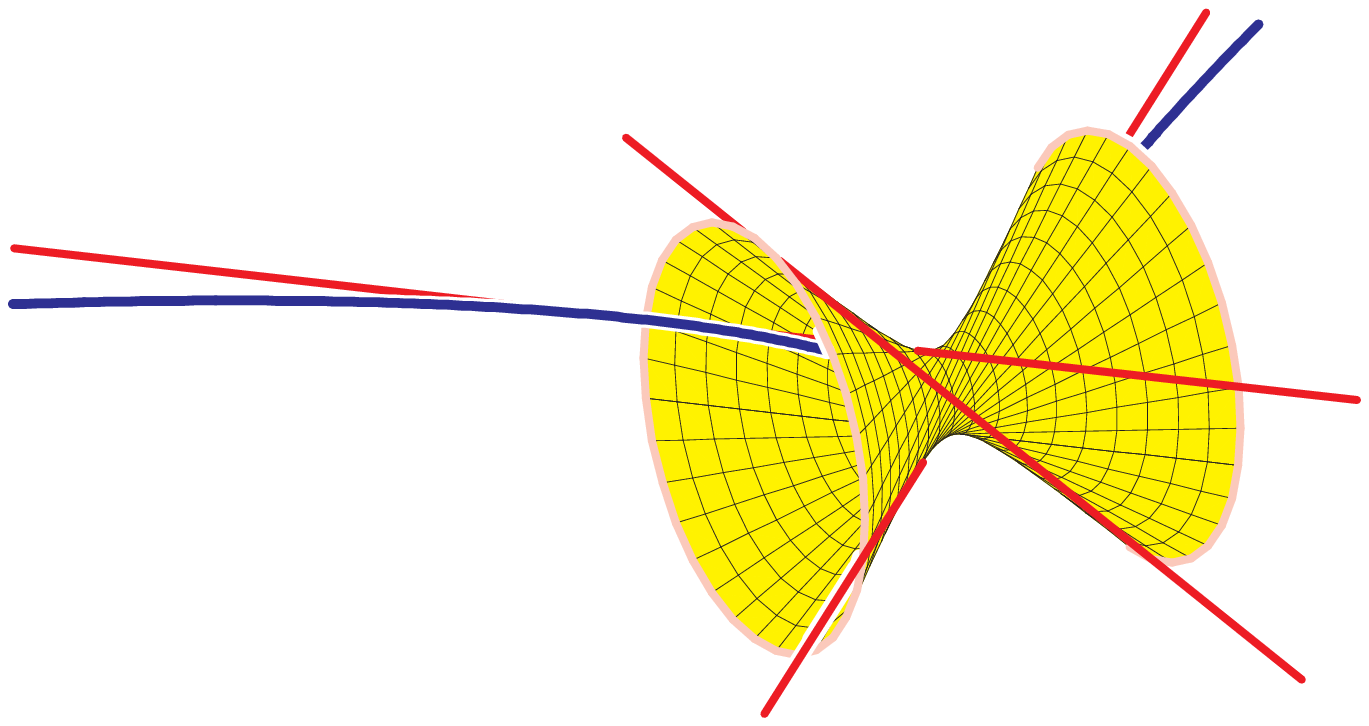}}
   \put(141.6,0){\Red{$\ell(t)$}} \put(295.2,0){\Red{$\ell(u)$}}
   \put(0,115.2){\Red{$\ell(s)$}} \put(0,74.9){\Blue{$\gamma$}}
   \put(270,115.2){\Brown{$Q$}}
  \end{picture}
 \]
 \caption{Quadric containing three lines tangent to the rational normal
 curve\label{F:3TanQuad}.} 
 \end{figure}

\begin{figure}[htb]
 \[
  \begin{picture}(240,102)(0,9.6)
   \put(0,0){\includegraphics[height=120pt]{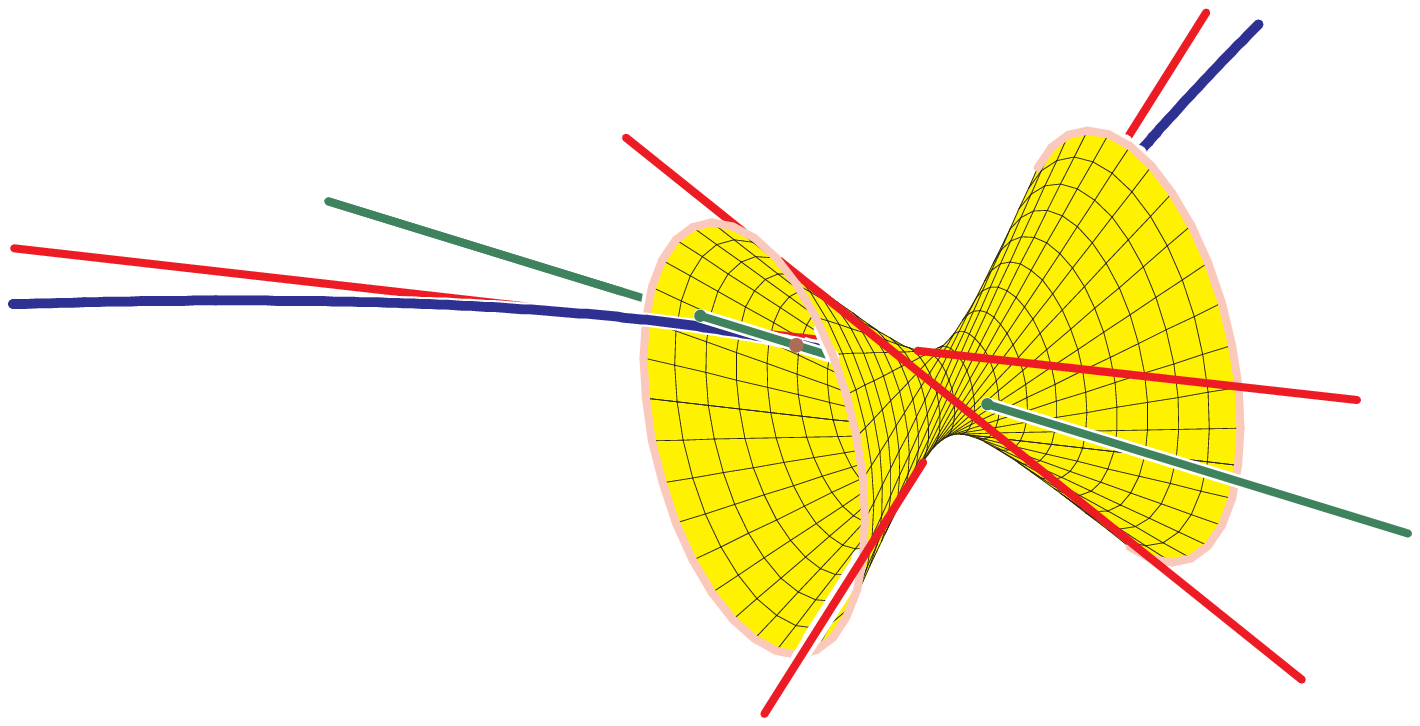}}
   \put(0,86.4){\Red{$\ell(s)$}} \put(105,0){\Red{$\ell(t)$}} 
   \put(194,0){\Red{$\ell(u)$}} 
   \put(0,54){\Blue{$\gamma$}} \put(36,93.6){\ForestGreen{$\lambda(v,w)$}}
  \end{picture}\quad
  \begin{picture}(198,110)(-24,-8)
   \put(0,0){\includegraphics[height=90pt]{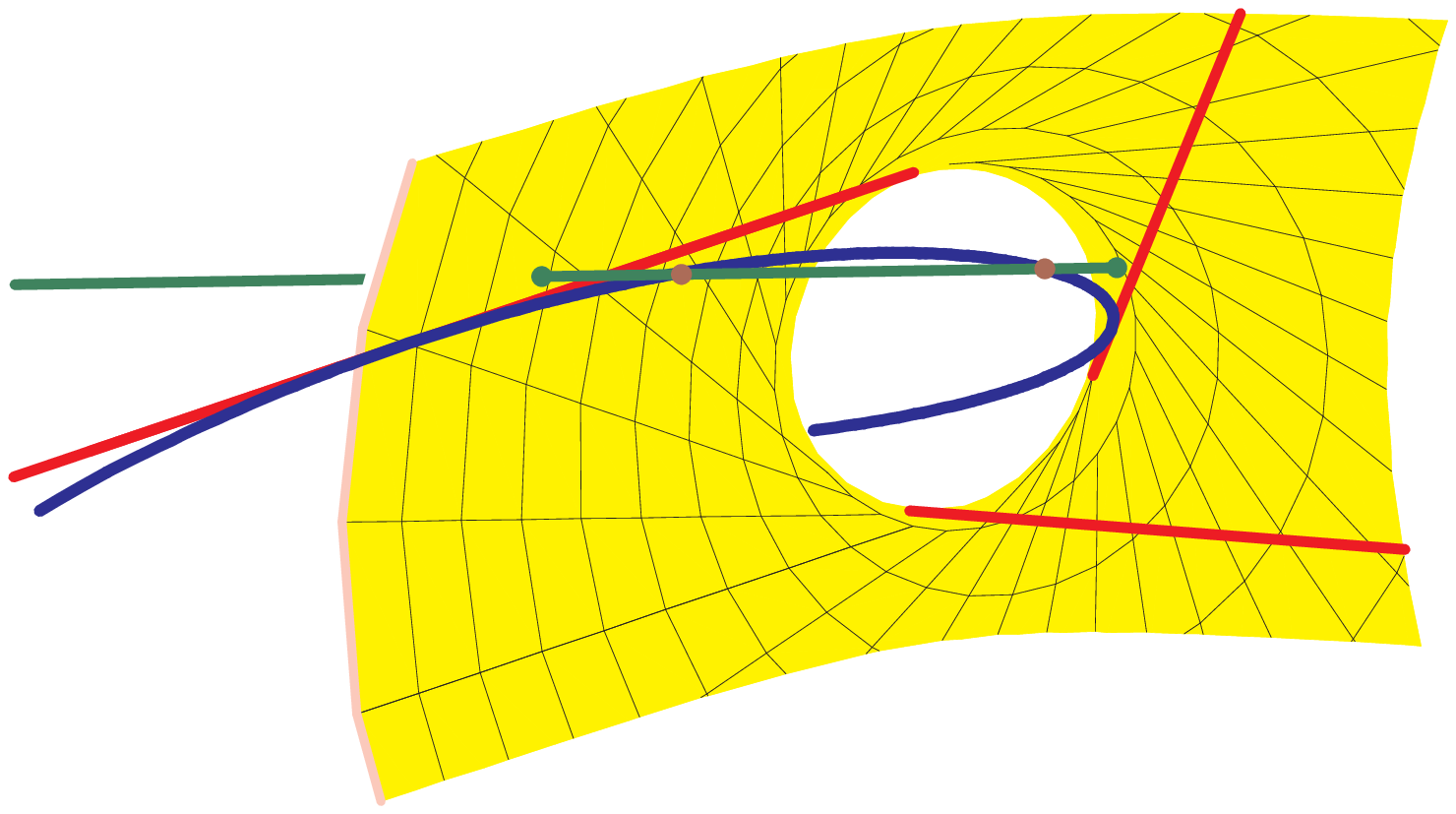} }
   \put(-14,44){\Red{$\ell(s)$}} \put(164.4,24){\Red{$\ell(t)$}} 
   \put(126,97){\Red{$\ell(u)$}}
   \put(5,24){\Blue{$\gamma$}} \put(-14,67){\ForestGreen{$\lambda(v,w)$}}
  \thicklines
    \put( 65,-10){\Brown{$\gamma(v)$}}  \put( 76,3){\Brown{\vector(0,1){54}}}
    \put(105,-10){\Brown{$\gamma(w)$}}  \put(117,3){\Brown{\vector(0,1){54}}} 
  \end{picture}
\] \caption{Two views of a secant line meeting $Q$\label{F:NC}.}
 \end{figure}

\begin{figure}[htb]
 \[
  \begin{picture}(360,170)(0,15)
   \put(0,0){\includegraphics[height=180pt]{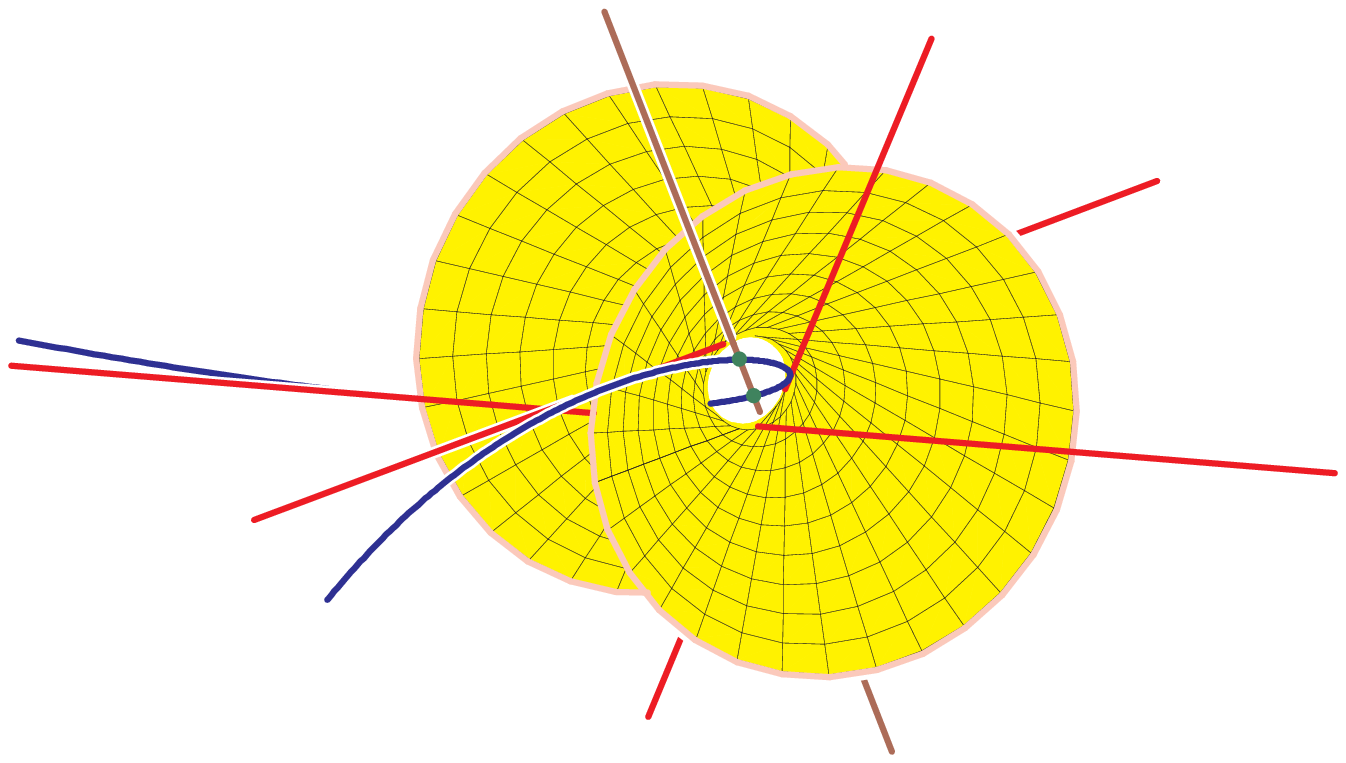} }
   \put(37,51){\Red{$\ell(s)$}} \put(300,50){\Red{$\ell(t)$}} 
   \put(228,168){\Red{$\ell(u)$}}
   \put(74.5,25,5){\Blue{$\gamma$}}\put(102,172.5){\Brown{$\lambda(v,w)$}}
   \thicklines
   \put(189,179){\ForestGreen{$\gamma(v)$}}  
       \put(198,174){\ForestGreen{\vector(-1,-4){18.3}}}
   \put(118,4){\ForestGreen{$\gamma(w)$}}  
       \put(133,15){\ForestGreen{\vector(2,3){46}}} 
  \end{picture}
 \]
 \caption{A secant line not meeting $Q$\label{F:CR}.}
 \end{figure}

The set of lines which meet the three tangent lines 
$\ell(s)$, $\ell(t)$, and $\ell(u)$ forms one ruling of a quadric surface $Q$
in $\P^3$. 
We display a picture of $Q$ and the ruling in
Figure~\ref{F:3TanQuad}, as well as the rational normal curve \Blue{$\gamma$}
with its three tangent lines.  
This is for a particular choice of $s$, $t$, and $u$, which is described below.
The lines meeting $\ell(s)$, $\ell(t)$, $\ell(u)$, and the secant line
$\lambda(v,w)$ correspond to the points where $\lambda(v,w)$ meets the quadric $Q$.
In Figure~\ref{F:NC}, we display a secant line $\lambda(v,w)$ which meets the
hyperboloid in two points, and therefore these choices for $v$ and $w$ give two real
flags in the intersection~\eqref{F234Int}.
There is also a secant line which meets the hyperboloid in no real points, and
hence in two complex conjugate points. 
For this secant line, both flags in the intersection~\eqref{F234Int} are 
complex. 
We show this configuration in Figure~\ref{F:CR}.

To investigate this failure of the Shapiro conjecture, first note
that any two parametrizations of two rational normal curves are conjugate
under a projective transformation of $\P^3$.
Thus it will be no loss to assume that the curve $\gamma$ has the 
parametrization
\[
   \gamma\ :\ t\ \longmapsto\ 
     [2,\, 12t^2-2, \, 7t^3+3t,\, 3t-t^3]\,.
\]
Then the lines tangent to $\gamma$ at the points $(s,t,u)=(-1,0,1)$
lie on the hyperboloid
\[
   x_0^2 - x_1^2 + x_2^2 - x_3^2\ =\ 0\,.
\]
If we parametrize the secant line $\lambda(v,w)$ as
$(\frac{1}{2}+l)\gamma(v) + (\frac{1}{2}-l)\gamma(w)$ and then 
substitute this into the equation for the hyperboloid, we obtain a quadratic
polynomial in $l,v,w$.
Its discriminant with respect to $l$ is
 \begin{equation}\label{Eq:F234_disc}
   16  (v - w)^2\, (2 v w + v + w) (3 v w + 1) (1 - v w) (v + w - 2 v w)\,.
 \end{equation}
We plot its zero-set in the square
$v,w\in[-2,2]$, shading the regions where the discriminant is negative. 
The vertical broken lines are $v,w=\pm1$, the diagonal line is $v=w$, 
the cross is the value of $(v,w)$ in Figure~\ref{F:NC}, and the dot is the value
in Figure~\ref{F:CR}.
 \begin{figure}[htb]
  \[
   \begin{picture}(162,154)
    \put(-2,-3){\includegraphics[height=160pt]{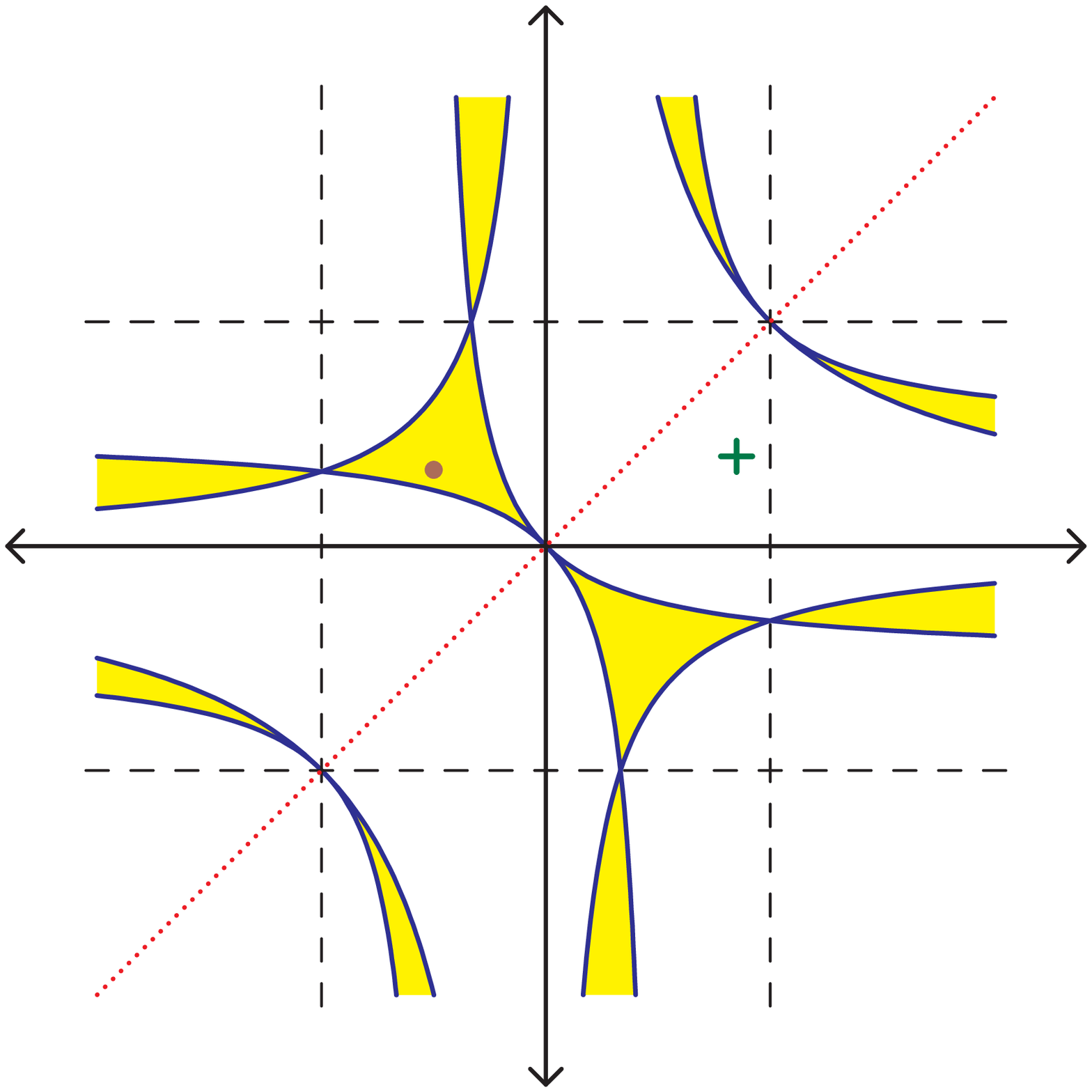}}
    \put(146, 83){$v$}   \put(83,146){$w$}
   \end{picture}
  \]
 \caption{Discriminant of the Schubert problem~\ref{F234Int}.\label{F:Discr}}
 \end{figure}
Observe that the discriminant is nonnegative if $(v,w)$ lies in one of the
squares $(-1,0)^2$, $(0,1)^2$, or if $(\frac{1}{v},\frac{1}{w})\in(-1,1)^2$
and it is positive in the triangles into which the line $v=w$ subdivides these
squares.  
Since $(s,t,u)=(-1,0,1)$, these squares are the values of $v$ and $w$ when
both lie entirely within one of the three intervals of $\RP^1$ determined by
$s,t,u$. 
If we allow M\"obius transformations of $\RP^1$, we deduce the
following proposition.

\begin{prop}\label{P:counter}
  The intersection~\eqref{F234Int} is transverse and consists only of real
  points if there are disjoint intervals $I_2$ and $I_3$ of\/ $\RP^1$ so that
  $s,t,u\in I_2$ and $v,w\in I_3$.
\end{prop}

While this example shows that the Shapiro conjecture is false, 
Proposition~\ref{P:counter} suggests that a refinement to the Shapiro conjecture
may hold. 
We will describe such a refinement and present
experimental evidence supporting it.

\section{Results}

Experimentation designed to test hypotheses is a
primary means of inquiry in the natural sciences.
In mathematics we use proof and example as our primary
means of inquiry.
Many mathematicians (including the authors) feel that they are 
striving to understand the nature of objects that inhabit a very
real mathematical reality.
For us, experimentation plays an important role in helping to  
formulate reasonable conjectures, which are then
studied and perhaps eventually decided.

We first discuss the conjectures which were informed by our 
experimentation that we describe in Section~\ref{S:Methods}.
Then we discuss the proof of these conjectures for the flag manifolds
$\Fl(n{-}2,n{-}1; n)$ by Eremenko, Gabrielov, Shapiro, and Vainshtein~\cite{EGSV}, and
an extension of our monotone conjecture which is suggested by their work.
Lastly, we present some examples from this experimentation which exhibit new
and interesting phenomena.

\subsection{Conjectures}\label{S:conjectures}

Let $\alpha=\{\alpha_1<\dotsb<\alpha_k\}$ and $n$ be positive integers with
$\alpha_k<n$. 
Recall that a permutation $w\in W^{\alpha}$ is Grassmannian if it has a single
descent, say at position $\alpha_l$.
Then the Schubert variety $X_wF_\bullet$ of $\Flan$ is the inverse image of the
Schubert variety $X_wF_\bullet$ of the Grassmannian $\Gr(\alpha_l,n)$.
Write $\delta(w)$ for the unique descent of a Grassmannian permutation $w$.

A Schubert problem  $(w_1,\dotsc,w_m)$ for $\Flan$ is 
\Blue{{\it Grassmannian}} if each permutation  $w_i$ is Grassmannian.
A list of points $t_1,\dotsc,t_m\in\RP^1$ is \Blue{{\it monotone}}
with respect to a Grassmannian Schubert problem $(w_1,\dotsc,w_m)$ if
the function 
\[
   t_i\ \longmapsto\ \delta(w_i)\ \in\{\alpha_1,\alpha_2,\dotsc,\alpha_k\}
\]
is monotone, when the ordering of the
$t_i$ is consistent with an orientation of $\RP^1$.
We also say that the ordered $m$-tuple $(t_1,\dotsc,t_m)$ is a monotone point of
$(\R\P^1)^m$.

This definition is invariant under the automorphism
group of $\R\P^1$, which consists of the real M\"obius
transformations and acts transitively on triples of
points on $\R\P^1$.
Viewing $\C^n$ as the linear space of homogeneous forms on $\P^1$ of degree
$n{-}1$ shows that an automorphism $\varphi$ of $\P^1$ induces a corresponding
automorphism $\varphi$ of $\C^n$ such that
$\varphi(\gamma(t))=\gamma(\varphi(t))$, and thus
$\varphi(F_\bullet(t))=F_\bullet(\varphi(t))$.
The corresponding automorphism $\varphi$ of $\Flan$ satisfies
$\varphi(X_w(t))=X_w(\varphi(t))$.
This was used in the discussion of Section~\ref{S:F124}.

\begin{rem}
 Conjecture~A of the Introduction involves a monotone choice of
 points for the Grassmannian Schubert problem
 $\bigl( \sigma_2^4,\, \sigma_3^4\bigr)$ on the 
 flag manifold $\mathbb{F}\ell(2,3;5)$.
 Indeed, $\calM_{\iota}$ is the set of matrices of the form
 \[
   \left[ \begin{matrix}
    1 & 0 & x_1 & x_2 & x_3\\
    0 & 1 & x_4 & x_5 & x_6\\
    0 & 0 &  1  & x_7 & x_8
    \end{matrix}\right]\ .
\]
 The equation $f(s;x)=0$ is the condition that $E_2(x)$ meets
 $F_3(s)$ non-trivially, and defines the Schubert variety $X_{\sigma_2}(s)$.
 Similarly, $g(s;x)=0$ defines the Schubert variety $X_{\sigma_3}(s)$. 
 The list of points at which $f$ and $g$ were evaluated in
 Conjecture~A is monotone.
\end{rem}

\begin{conj}\label{C:Main}
 Suppose that $(w_1,\dotsc,w_m)$ is a Grassmannian Schubert problem for $\Flan$.
 Then the intersection 
\begin{equation}\label{AA}
     X_{w_1}(t_1)\cap X_{w_2}(t_2)\cap 
    \dotsb\cap X_{w_m}(t_m)\,,
\end{equation}
 is transverse with all points of intersection real, if the points
 $t_1,\dotsc,t_m\in\RP^1$ are monotone with respect to $(w_1,\dotsc,w_m)$.
\end{conj}

We make a weaker conjecture which drops the claim of transversality.

\begin{conj}\label{C:Main-NT}
 Suppose that $(w_1,\dotsc,w_m)$ is a Grassmannian Schubert problem for $\Flan$.
 Then the intersection~\eqref{AA} has all points real, if the points
 $t_1,\dotsc,t_m\in\RP^1$ are monotone with respect to $(w_1,\dotsc,w_m)$.
\end{conj}

\begin{rem}
 The example of Section~\ref{S:F124} illustrates both Conjecture~\ref{C:Main}
 and its limitation.
 The condition on disjoint intervals $I_2$ and $I_3$ of Proposition~\ref{P:counter} is
 equivalent to the pointss being monotone.
 The shaded regions in Figure~\ref{F:Discr}, which are the points that give no real
 solutions, contain no monotone lists of points.
\end{rem}

If $\Flan$ is a Grassmannian, then every choice of points is monotone, so
Conjecture~\ref{C:Main} includes the Shapiro conjecture for
Grassmannains as a special case. 
Our experimentation systematically investigated the original Shapiro conjecture
for flag manifolds, with a focus on this monotone conjecture.
We examined 590 such Grassmannian Schubert problems on 29 different flag
manifolds.
In all, we verified that each of more than 158 million specific monotone
intersections of the form~\eqref{AA} had all solutions real.
We find this to be overwhelming evidence in support of our monotone conjecture.

Indeed, the set of points
$(t_1,\dotsc,t_m)\in(\P^1)^m$ where the intersection~\eqref{AA} is not
transverse is the \Blue{{\it discriminant}}
\Blue{$\Sigma$} of the corresponding Schubert problem.
This is a hypersurface, unless the intersection
is never transverse.
The number of real solutions is constant on each connected component of the
complement of the discriminant.
Conjecture~\ref{C:Main} asserts that the set of monotone points lies 
entirely within the region where all solutions are real.
Our computations show that the discriminant is a hypersurface for the
Grassmannian Schubert problems we considered, and none of the 158 million
monotone points we considered was contained in a non-maximal
component in which not all solutions were real.
While this does not prove Conjecture~\ref{C:Main} for these problems, it
places severe restrictions on the location of the non-maximal components of
the complement of the discriminant.

For a given flag manifold, it suffices to know Conjecture~\ref{C:Main-NT} for 
\Blue{{\it simple Schubert problems}}, which involve only simple (codimension 1)
Schubert conditions. 
As simple Schubert conditions are Grassmannian, 
Conjectures~\ref{C:Main} and~\ref{C:Main-NT} apply to simple Schubert problems.

\begin{thm}\label{T:simple}
 Suppose that Conjecture~$\ref{C:Main-NT}$ holds for all simple Schubert
 problems on a given flag manifold $\Flan$.
 Then Conjecture~$\ref{C:Main-NT}$ holds for all Grassmannian Schubert
 problems on any flag manifold $\Fl(\beta;n)$ where $\beta$ is a subsequence of
 $\alpha$. 
\end{thm}

We prove Theorem~\ref{T:simple} when $\beta=\alpha$ in
Section~\ref{S:T_simple} and the general case in
Section~\ref{S:projections}. 

We give two further and successively stronger conjectures which are supported by our
experimental investigation.
The first ignores the issue of reality and concentrates only on the 
transversality of an intersection.

\begin{conj}\label{C:transverse}
 If $(w_1,\dotsc,w_m)$ is a Grassmannian Schubert problem for $\Flan$ and
 the points $t_1,\dotsc,t_m\in\RP^1$ are monotone with respect to
 $(w_1,\dotsc,w_m)$, then the intersection~\eqref{AA} is transverse.
\end{conj}

Since the set of monotone points is connected, Conjecture~\ref{C:transverse}
asserts that it lies 
in a single component of the complement of the discriminant.
Since a main result of~\cite{So00b} is that Conjecture~\ref{C:Main} holds
for simple Schubert problems when the points $t_1,\dotsc,t_m$ are
sufficiently clustered together, Conjecture~\ref{C:transverse} implies
Conjecture~\ref{C:Main}, for simple Schubert problems. 
Then Theorem~\ref{T:simple} implies Conjecture~\ref{C:Main-NT}, and the
transversality assertion of Conjecture~\ref{C:transverse} implies 
Conjecture~\ref{C:Main}, without any restriction on the Grassmannian Schubert
problem. 

\begin{thm}\label{T:Implies}
 Conjecture~$\ref{C:transverse}$ implies Conjecture~$\ref{C:Main}$.
\end{thm}

Conjecture~\ref{C:transverse} states that for a
Grassmannian Schubert problem $w$, the discriminant $\Sigma$ contains no points
$(t_1,\dotsc,t_m)$ that are monotone with respect to $w$. 
In our experimentation, we kept track of the 
non-transverse intersections.
None came from monotone points for a
Grassmannian Schubert problem.
In contrast, there were several hunderd such non-transverse intersections
encountered involving non-monotone choices of points.
While this does not rule out the existence of monotone choices of points giving
a non-transverse intersection, it does suggest that it is highly unlikely.

In every case that we have computed, the discriminant is defined by a
polynomial having a special form which shows that $\Sigma$ contains no points
that are monotone with respect to $w$.
We explain this.
The set $\Sigma\cap \R^m$ is defined by a single 
\Blue{{\it discriminant polynomial}} $\Delta_w(t_1,\dotsc,t_m)$, that is
well-defined up to multiplication by a scalar.
The set of monotone points $(t_1,\dotsc,t_m)\in\R^m$ with
respect to $w$ has many components.
Consider the union of components defined
by the inequalities 
 \begin{equation}\label{E:monotone}
   t_i\neq t_j\quad\mbox{if}\quad i\neq j
   \qquad\mbox{and}\qquad
   t_i<t_j\quad\mbox{whenever}\quad
   \delta(w_i)<\delta(w_j)\,.
 \end{equation}

For the example of Section~\ref{S:F124}, the region of monotone points 
is where $v,w$ lie in one of the three intervals of $\R\P^1$ defined by $s,t,u$.
As we argued there, we may assume that $(s,t,u)=(-1,0,1)$ and so 
$v,w$ must lie in one of the three disjoint intervals $(-1,0)$, $(0,1)$, or $(1,-1)$ on
$\R\P^1$, where the last interval contains $\infty$.
Since any one of these intervals is transformed into any other by a M\"obius
transformation, it suffices to consider the interval $(0,1)$, which is defined
by the inequalities
\[
   0\ <\ v\,, w\,,\qquad\mbox{and}\qquad 0\ <\ 1-v\,,\, 1-w\,.
\]
Note that 
\begin{eqnarray*}
   1-vw  &=&  1{-}w + w(1{-}v)    \\
  v+w-2vw&=& v(1{-}w) +  w(1{-}v) \,, 
\end{eqnarray*}
which shows that the discriminant~\eqref{Eq:F234_disc} is positive if $v\neq w$
and $0<v,w<1$.

We conjecture that the discriminant always has such a form for which its
positivity (or negativity) on the set~\eqref{E:monotone} of monotone points
is obvious. 
More precisely, 
suppose that $w=(w_1,\dotsc,w_m)$ is a Grassmannian Schubert problem for $\Flan$.
Set 
\[
   S\ :=\ \{ t_i-t_j \mid \delta(w_i)>\delta(w_j)\}\,.
\]
Then the set~\eqref{E:monotone} of monotone points is
\[
   \{ t=(t_1,\dotsc,t_m)\mid g(t)\geq 0\mbox{\quad for\quad} g\in S\}\,.
\]
Writing $S=\{g_1,\dotsc,g_l\}$, the \Blue{{\it preorder}}
generated by $S$ is the set of polynomials of the form 
\[
   \sum_\varepsilon  c_\varepsilon g_1^{\varepsilon_1} 
          g_2^{\varepsilon_2}\dotsc  g_l^{\varepsilon_l} \,,
\]
where each $\varepsilon_i\in\{0,1\}$ and each coefficient $c_\varepsilon$ is a sum of
squares of polynomials.
Every polynomial in the preorder generated by $S$ is obviously
positive on the set~\eqref{E:monotone} of monotone points, but not every
polynomial that is positive on that set lies in the
preorder, at least when $m\geq 5$.
Indeed, suppose that $\delta(w_1)\leq\delta(w_2)\leq\dotsb\leq\delta(w_m)$.
Using the automorphism group of $\RP^1$, we may assume that 
$t_1=\infty$, $t_2=-1$, $t_3=0$.
Then the set~\eqref{E:monotone} are those 
$(t_4,\dotsc,t_m)$ such that $0<t_4<\dotsb<t_m$.
This contains a 2-dimensional cone when $m\geq 5$, so the preorder of 
polynomials which are positive on this set is not a finitely
generated preorder~\cite[\S6.7]{Sch00}.

\begin{conj}\label{C:Discriminant}
 Suppose that $(w_1,\dotsc,w_m)$ is a Grassmannian Schubert problem for $\Flan$.
 Then its discriminant $\Delta_w$ (or its negative) lies in the preorder 
 generated by the polynomials 
\[
   S\ :=\ \{ t_i-t_j \mid \delta(w_i)>\delta(w_j)\}\,.
\]
\end{conj}

We showed that this holds for the problem of Section~\ref{S:F124}.
Conjecture~\ref{C:Discriminant} generalizes a conjecture made in~\cite{So00a}
that the discriminants for Grassmannians are sums of squares.

Since Conjecture~\ref{C:Discriminant} implies that the discriminant is
nonvanishing on monotone choices of points, it implies 
Conjecture~\ref{C:transverse}, and so by Theorem~\ref{T:Implies}, it implies the
original Conjecture~\ref{C:Main}.
We record this fact.

\begin{thm}\label{T:ImpliesII}
 Conjecture~$\ref{C:Discriminant}$ implies Conjecture~$\ref{C:Main}$.
\end{thm} 

We give some additional evidence in favor of Conjecture~\ref{C:Discriminant} in
Section~\ref{S:discriminants}.

\subsection{The result of Eremenko, Gabrielov, Shapiro, and
  Vainshtein}

Conjecture \ref{C:Main} for
${\mathbb F}\ell(n{-}2,n{-}1;n)$ follows from a result of 
Eremenko {\it et.~al}~\cite{EGSV}.  
We discuss this for simple Schubert problems, from which the general
case follows, by Theorem~\ref{T:simple}.  

There are two types of simple Schubert 
varieties in ${\mathbb F}\ell(n{-}2,n{-}1;n)$,
 \begin{eqnarray*}
   X_{\sigma_{n{-}2}}F_\bullet
    & := & \{(E_{n{-}2}\subset E_{n{-}1})\mid E_{n{-}2}\cap F_2\neq \{0\}\}\,,
    \makebox[1pt][l]{\qquad and}\\\ 
   X_{\sigma_{n{-}1}}F_\bullet
    & := & \{(E_{n{-}2}\subset E_{n{-}1})\mid E_{n{-}1}\supset F_1\}\,.
 \end{eqnarray*}
When $n=4$, these are the Schubert varieties $X_{\sigma_2}F_\bullet$ and
$X_{\sigma_3}F_\bullet$ of Section~\ref{S:F124}. 

Consider the Schubert intersection
 \begin{equation}\label{codim1}
  X_{\sigma_{n{-}2}}(t_1)\cap\dotsb\cap X_{\sigma_{n{-}2}}(t_p)~\cap~
  X_{\sigma_{n{-}1}}(s_1)\cap\dotsb\cap X_{\sigma_{n{-}1}}(s_q)
 \end{equation}
where $t_1,\dotsc,t_p$ and $s_1,
\dotsc,s_q$ are distinct points in $\RP^1$ and $p+q=2n-1$ with
$0<q\leq n$.
As in Section~\ref{S:F124}, this Schubert problem is equivalent to one on the 
Grassmanian $\Gr(n{-}2,n)$ of codimension 2 planes.  
The condition that $E_{n{-}1}$ contains 
each of the 1-dimensional linear subspaces ${\rm span}\{\gamma(s_i)\}$ 
for $i=1,\dotsc,q$ implies that $E_{n{-}1}$
contains the secant plane 
$W={\rm span}\{\gamma(s_i)|i=1,\dotsc,q\}$ of dimension $q$.
This forces the condition that $\dim W\cap E_{n{-}2}\geq q{-}1$,
so that $E_\bullet\in X_{\tau}W$, where $\tau$ is the Grassmannian permutation
\[
   (1,\ 2,\ \dotsc,\ n{-}q,\ \ \ \ n{-}q{+}2,\ \dotsc,\ 
    n{-}1,\ \ n{-}q{+}1,\ \ n)\ .
\]

One the other hand, when $\dim W\cap E_{n{-}2}=q{-}1$, we can recover the 
hyperplane $E_{n{-}1}$ by setting $E_{n{-}1}:=W+E_{n{-}2}$.  
Thus the Schubert problem~\eqref{codim1}
reduces to a Schubert problem on $\Gr(n{-}2,n)$ of the form
 \begin{equation}\label{secSchInt}
   X_{\sigma_{n{-}2}}(t_1)\cap\dotsb\cap
   X_{\sigma_{n{-}2}}(t_p)~\cap~X_{\tau}W\,.
 \end{equation}
Using the results of~\cite{EG02a}, Eremenko, Gabrielov, Shapiro and Vainshtein
show that the intersection~\eqref{secSchInt} has only real points, when the
given points $t_1,\dotsc,t_p,s_1,\dotsc,s_q$ are monotone with respect to the
Schubert problem $(\sigma_{n-2}^p,\,\sigma_{n-1}^q)$.

This suggests a generalization
of Conjecture~\ref{C:Main} to flags of subspaces which
are secant to the rational normal curve $\gamma$. 
Let $S:=(s_1,s_2,\dotsc,s_n)$ be $n$ distinct points in $\P^1$ and for each
$i=1,\dotsc,n$, let $F_i(S):={\rm span} \{\gamma(s_1),\dotsc,\gamma(s_i)\}$.
These subspaces form the flag $F_\bullet(S)$ which is 
\Blue{{\it secant to $\gamma$ at $S$}}.
A list $(S_1,\dotsc,S_m)$, of sets of $n$ distinct points in $\RP^1$ is
\Blue{{\it monotone}} with respect to a Grassmannian Schubert problem
$(w_1,\dotsc,w_m)$ if 
 \begin{enumerate}
  \item There exists a collection of disjoint intervals 
        $I_1,\dotsc,I_m$ of $\RP^1$ with 
        $S_i\subset I_i$ for each $i=1,\dotsc,m$, and 
  \item If we choose points  $t_i\in I_i$ for $i=1,\dotsc,m$, then 
        $(t_1,\dotsc,t_m)$ is monotone with respect
        to the Grassmannian Schubert problem $w$.
        This notion does not depend upon the choice of points, as the
        intervals are disjoint.
\end{enumerate}

\begin{conj}\label{C:secant}
 Given a Grassmannian Schubert problem $(w_1,\dotsc,w_m)$ for $\Flan$, 
 the Schubert intersection 
\[
     X_{w_1}F_\bullet(S_1)\cap X_{w_2}F_\bullet(S_2)\cap 
    \dotsb\cap X_{w_m}F_\bullet(S_m)\,,
\]
 is transverse with all points of intersection real, if the list of subsets
 $(S_1,\dotsc,S_m)$ of $\RP^1$ is monotone with respect to $(w_1,\dotsc,w_m)$.
\end{conj}

Conjecture~\ref{C:secant} was formulated in the case when the flag manifolds are
Grassmannians in~\cite{EGSV}, where monotonicity was called well-separatedness. 
The main result in that paper is its proof for the Grassmannian $\Gr(n{-}2,n)$.
A collection $U_1,\dotsc,U_r$ of subsets of $\RP^1$ is \Blue{{\it well-separated}} if
there are disjoint intervals $I_1,\dotsc,I_r$ of $\RP^1$ with $U_i\subset I_i$ for 
$i=1,\dotsc,r$.

\begin{prop}[Eremenko,~{\it et.~al}~\protect{\cite[Theorem~1]{EGSV}}]
\label{P:secant}
 Suppose that $U_1,\dotsc,U_r$ is a well-separated collection of finite subsets
 of $\RP^1$ consisting of $2n-2+r$ points, and with no $U_i$ consisting of a
 single point.
 Then there are finitely many codimension $2$ planes meeting each of the planes
 $\mbox{span}\{\gamma(U_i)\}$ for $i=1,\dotsc,r$, and all are real.
\end{prop}

The numerical condition that there are $2n-2+r$ points and that no $U_i$ is a
singleton ensures that there will be finitely many codimension 2 planes meeting
the subspaces $\mbox{span}\{\gamma(U_i)\}$.
To see how this implies that the intersections~\eqref{secSchInt}
and~\eqref{codim1} consist only of real points, let $r=p+1$ and set
$U_j:=\{t_j,u_j\}$, where the point $u_j$ is close to the point $t_j$ for
$j=1,\dotsc,p$ and also set $U_{p+1}:=\{s_1,\dotsc,s_q\}$. 
For each $j=1,\dotsc,p$, the limit
\[
 \lim_{u_j\rightarrow t_j}\mbox{span}\{\gamma(U_j)\}
\]
is the 2-plane osculating the rational normal curve at $t_j$. 
The condition that the subsets $U_1,\dotsc,U_{p{+}1}$ are are well-separated
implies that the points $\{s_1,\dotsc,s_q,t_1,\dotsc,t_p\}$
are monotone with respect to the Schubert problem 
$(\sigma_{n{-}2}^p,\sigma_{n{-}1}^q)$.
Thus the intersection~\eqref{secSchInt} is a limit of intersections
of the form in Proposition~\ref{P:secant}, and hence consists only of real points.
This gives the following corollary to Proposition~\ref{P:secant},
also proven in~\cite{EGSV}.

\begin{cor}
 Suppose that there exist disjoint intervals 
 $I\supset\{s_1,\dotsc,s_q\}$ and $J\supset\{t_1,\dotsc,t_p\}$.
 Then all codimension $2$ planes in the intersection~$\eqref{codim1}$ are real.  
 Thus all flags $E_\bullet\in{\mathbb F}\ell(n{-}2 , n{-}1;n)$ in
 the intersection~$\eqref{secSchInt}$ are real.
\end{cor}

We have not yet investigated Conjecture~\ref{C:secant},
and the results of~\cite{EGSV} are the only evidence currently in its favor.
We believe that experimentation testing this conjecture, in the spirit
of the experimentation described in Section 4, is a natural and
worthwhile next step.

\subsection{Examples}
While the original goal of our experimentation was to study
Conjecture~\ref{C:Main}, 
this project became a general study of Schubert
intersection problems on small flag manifolds.  
Here, we report on some new and interesting phenomena which we
observed, beyond support for Conjecture~\ref{C:Main}.

We first discuss some of the Schubert problems that we investigated, presenting in
tabular form the data from our experimentation on those problems.
Some of these appear to present new or interesting phenomena beyond
Conjecture~\ref{C:Main}.
We next discuss some phenomena that we observed in our data, and which we can 
establish rigorously.
One is the smallest enumerative problem that we know of with an
unexpectedly small Galois group~\cite{Ha79,Va04}, and the other is a Schubert
problem for which the intersection is not transverse, when the given flags
osculate the rational normal curve.\medskip 

A Schubert intersection of the form
\[
  X_{w_1}(t_1)\cap X_{w_2}(t_2)\cap\dotsb\cap X_{w_m}(t_m)
\]
may be encoded by labeling each point $t_i\in\RP^1$ with the corresponding Schubert
condition $w_i$.
The automorphism group of $\RP^1$ acts on the flag variety $\Flan$, and hence on
collections of labeled points.
A coarser equivalence which captures the combinatorics of the arrangement of Schubert
conditions along $\RP^1$ is isotopy, and isotopy classes of such labeled points are
called \Blue{{\it necklaces}}, which are the different arrangements of
$m$ beads labeled with $w_1,\dotsc,w_m$ and strung on the circle $\RP^1$.
Our experimentation was designed to study how the number of real
solutions to a Schubert problem was affected by the necklace.
\Blue{{\it Monotone necklaces}} are necklaces corresponding to monotone choices
of points.

To that end, we kept track of the number of real solutions to a Schubert problem
by the associated necklace, and have archived the results in linked web pages
available at {\tt www.math.tamu.edu/\~{}sottile/pages/Flags/}.
Section~\ref{S:Methods} discusses how these computations were carried out.
While Conjecture~\ref{C:Main} is the most basic assertion that we believe is
true, there were many other phenomena, both general and specific, that our
experimentation uncovered.
We describe some of them below.
Conjecture~\ref{C:New_Main} and Theorem~\ref{T:Projections} are some others. 
Our data contain many more interesting examples, and
invite the interested reader browse the data online.

\subsubsection{Conjecture~$\ref{C:Main}$}
Table~\ref{table:12-flag} shows the data from computing 3.2 million 
instances of the Schubert problem $(\Blue{\sigma_2}^4,\Magenta{\sigma_3}^4)$ on 
$\Fl(2,3;5)$ underlying Conjecture~A from the
 \begin{table}[htb]
  \begin{tabular}
   {|c||c|c|c|c|c|c|c|}\hline
   {Necklace} & \multicolumn{7}{c|}{Number of Real Solutions\rule{0pt}{13pt}}\\
   \cline{2-8}
        &0&2&4&6&8&10&12\rule{0pt}{13pt}\\\hline\hline
   \Blue{2222}\Magenta{3333}
     &0 & 0 & 0 & 0 & 0 & 0 & 400000\\\hline
   \Blue{22}\Magenta{3}\Blue{22}\Magenta{333}   
     &0 & 0 & 118 & 65425 & 132241 & 117504 & 84712\\ \hline
   \Blue{222}\Magenta{33}\Blue{2}\Magenta{33}
     &0 & 0 & 104 & 65461 & 134417 & 117535 & 82483 \\ \hline
   \Blue{22}\Magenta{33}\Blue{22}\Magenta{33}
     &0 & 0 & 1618 & 57236 & 188393 & 92580 & 60173  \\ \hline
   \Blue{22}\Magenta{3}\Blue{2}\Magenta{33}\Blue{2}\Magenta{3}
     &0 & 0 & 25398 & 90784 & 143394 & 107108 & 33316 \\ \hline
   \Blue{22}\Magenta{33}\Blue{2}\Magenta{3}\Blue{2}\Magenta{3}
     &0 & 2085 & 79317 & 111448 & 121589 & 60333 & 25228\\ \hline
   \Blue{222}\Magenta{3}\Blue{2}\Magenta{333}
     &0 & 7818 & 34389 & 58098 & 101334 & 81724 & 116637 \\ \hline
   \Blue{2}\Magenta{3}\Blue{2}\Magenta{3}\Blue{2}\Magenta{3}\Blue{2}\Magenta{3}
     & 15923 & 41929 & 131054 & 86894 & 81823 & 30578 & 11799 \\ \hline
  \end{tabular}\vspace{5pt}

  \caption{The Schubert problem $(\Blue{\sigma_2}^4,\Magenta{\sigma_3}^4))$ on
  $\Fl(2,3;5)$.}  
  \label{table:12-flag}
 \end{table}
%
%
Introduction. 
Each row corresponds to a necklace, and the entries record how often a
given number of real solutions was observed for the corresponding necklace.  
Representing the Schubert conditions $\Blue{\sigma_2}$ and $\Magenta{\sigma_3}$  by
their subscripts, we may write each necklace
linearly as a sequence of 2s and 3s. 
The only monotone necklace is in the first row, and Conjecture~\ref{C:Main}
predicts that any intersection with this necklace will have all 12 solutions
real, as we observe.

The other rows in this table are equally striking.
It appears that there is a unique necklace for which it is possible that no
solutions are real, and for five of the necklaces, the minimum number of
real solutions is 4. 
The rows in this and all other tables are ordered to highlight this feature.
Every row has  a non-zero entry in its last column.
This implies that for every necklace, there is a choice of points on $\RP^1$
with that necklace for which all 12 solutions are real.
Since this is a simple Schubert problem, that feature is a 
consequence of Corollary~2.2 of~\cite{So99a}.

Table~\ref{table:new12-flag} shows data from a related problem
$(\Purple{\sigma_1}^2,\Blue{\sigma_2}^3,
  \ForestGreen{\sigma_3}^3,\Red{\sigma_4}^2)$ with 12 solutions.
We only computed three necklaces for this problem, as it has 1,272 necklaces.
 \begin{table}[htb]
  \begin{tabular}
   {|c||c|c|c|c|c|c|c|}\hline
   {Necklace} & \multicolumn{7}{c|}{Number of Real Solutions\rule{0pt}{13pt}}\\
   \cline{2-8}
        &0&2&4&6&8&10&12\rule{0pt}{13pt}\\\hline\hline
%
  \Purple{11}\Blue{222}\ForestGreen{333}\Red{44} 
       & 0 &  0 &  0 &   0 &    0 &   0 & 10000 \\\hline
  \Purple{11}\Blue{222}\Red{44}\ForestGreen{333} 
       & 0 &  0 &  0 &   0 &    0 &   0 & 10000 \\\hline  
  \Purple{11}\ForestGreen{333}\Blue{222}\Red{44} 
       & 0 & 102 &462 & 1556 & 3821 & 2809 & 1250 \\\hline
  \end{tabular}\vspace{5pt}
  \caption{The Schubert problem $(\Purple{\sigma_1}^2,\Blue{\sigma_2}^3,
  \ForestGreen{\sigma_3}^3,\Red{\sigma_4}^2)$  on  $\Fl(1,2,3,4;5)$.}  
  \label{table:new12-flag}
 \end{table}
In the necklaces, $i$ represents the Schubert condition $\sigma_i$.
The only monotone necklace is in the first row.
While the second row is not monotone, it appears to have only real solutions.
A similar phenomenon (some non-monotone necklaces having only real solutions)
was observed in other Schubert problems involving 4- and 5-step flag manifolds.  
This can be seen in the example of Table~\ref{table:11352}, as well as
the third part of Theorem~\ref{T:discriminants}.

Table~\ref{table:11352} shows data from
the problem $(\Purple{\sigma_1}^2,\Blue{\sigma_2}^2, \ForestGreen{246},
 \ForestGreen{\sigma_3},\Red{\sigma_4}^2,\Magenta{\sigma_5}^2)$ on
$\Fl(1,2,3,4,5;6)$ with 8 solutions.
In the necklaces, $i$ represents $\sigma_i$ and $\ForestGreen{C}$ represents the
Grassmannian condition $\ForestGreen{246}$ with descent at 3.
 \begin{table}[htb]
  \begin{tabular}
   {|c||c|c|c|c|c|}\hline
   {Necklace} & \multicolumn{5}{c|}{Number of Real Solutions\rule{0pt}{13pt}}\\
   \cline{2-6}
        &0&2&4&6&8\rule{0pt}{13pt}\\\hline\hline
    \Purple{11}\Blue{22}\ForestGreen{C3}\Red{44}\Magenta{55}  
                    & 0 &    0 &     0 &     0 & 50000\\\hline
    \Purple{11}\ForestGreen{C3}\Red{44}\Magenta{55}\Blue{22}  
                    & 0 &    0 &     0 &     0 & 50000\\\hline
    \Purple{11}\Blue{22}\ForestGreen{C3}\Magenta{55}\Red{44}  
                    & 0 &    0 &     0 &     0 & 50000\\\hline
    \Purple{11}\ForestGreen{C3}\Magenta{55}\Red{44}\Blue{22}  
                    & 0 &    0 &     0 &     0 & 50000\\\hline
    \Purple{11}\Magenta{55}\Blue{22}\ForestGreen{C3}\Red{44}  
                    & 0 &    0 &     0 &  3406 & 46594\\\hline
    \Purple{11}\ForestGreen{C3}\Magenta{55}\Blue{22}\Red{44}  
                    & 0 &  0 &  5401 & 24714&  19885\\\hline
    \Purple{11}\Magenta{55}\ForestGreen{C3}\Red{44}\Blue{22}  
                    & 0 & 0 &  6347 & 19567 & 24086\\\hline
    \Purple{11}\Blue{22}\Magenta{55}\ForestGreen{C3}\Red{44}  
                    & 0 & 0 &  7732 & 23461&  18807\\\hline
    \Purple{11}\ForestGreen{C3}\Red{44}\Blue{22}\Magenta{55}  
                    & 0 & 0&  12437&  20396 & 17167\\\hline
    \Purple{11}\Red{44}\Blue{22}\ForestGreen{C3}\Magenta{55}  
                    & 0 & 0&  12508 & 19177 & 18315\\\hline
    \Purple{11}\Red{44}\Magenta{55}\Blue{22}\ForestGreen{C3}  
                    & 0 & 0 & 15109 & 25418 &  9473\\\hline
    \Purple{11}\Magenta{55}\Red{44}\Blue{22}\ForestGreen{C3}  
                    & 0 & 0 & 17152 & 23734 &  9114\\\hline
    \Purple{1}\ForestGreen{3}\Magenta{5}\Blue{2}\Red{4}\Purple{1}\ForestGreen{C}\Magenta{5}\Blue{2}\Red{4} 
               & 298 & 7095 & 18280 & 17871 & 6456\\\hline
  \end{tabular}\vspace{5pt}
  \caption{The Schubert problem $(\Purple{\sigma_1}^2,\Blue{\sigma_2}^2, \ForestGreen{246},
 \ForestGreen{\sigma_3},\Red{\sigma_4}^2,\Magenta{\sigma_5}^2)$  on  $\Fl(1,2,3,4,5;6)$.}  
  \label{table:11352}
 \end{table}
We only computed 13 necklaces for this problem, as it has 11,352 necklaces.
Note that three non-monotone necklaces have only real solutions, one has at
least 6 solutions, and 7 have at least 4 real solutions.

\subsubsection{Apparent lower bounds}
In the last section, we noted that the lower bound on the
number of real solutions seems to depend upon the necklace.
We also found many Schubert problems with an
apparent \Red{\emph{lower bound}} which holds for \Red{{\it all}} necklaces.  
For example, Table~\ref{table:lower} is for the  Schubert problem 
$(\Blue{\sigma_3}, (\Blue{1362})^2, \, 
  \Magenta{\sigma_4}^2,\Magenta{1346})$ on $\Fl(3,4;7)$,
which has degree 10. 
We only display 4 of the 16 necklaces for this problem.
%
 \begin{table}[htb]
  \begin{tabular}
  {|c||c|c|c|c|c|c|}\hline
   {Necklace} & \multicolumn{6}{c|}{Number of Real Solutions\rule{0pt}{13pt}}\\
   \cline{2-7}
        &0&2&4&6&8&10\rule{0pt}{13pt}\\\hline\hline
  \Blue{abb}\Magenta{ccd} 
           & 0 & 0 &   0   &   0   &    0  & 100000 \\\hline 
  \Blue{a}\Magenta{c}\Blue{bb}\Magenta{cd}
           & 0 & 0 &   0   & 16722 & 50766 &  32512 \\\hline 
  \Blue{a}\Magenta{cc}\Blue{bb}\Magenta{d}
           & 0 & 0 & 11979 & 26316 & 29683 &  32022 \\\hline 
  \Blue{a}\Magenta{c}\Blue{b}\Magenta{d}\Blue{b}\Magenta{c}
           & 0 & 0 & 27976 & 34559 & 26469 &  10996 \\\hline 
\end{tabular} \vspace{5pt}
  \caption{The Schubert problem $(\Blue{\sigma_3}, (\Blue{1362})^2, \, 
             \Magenta{\sigma_4}^2, \Magenta{1346})$ on
  $\Fl(3,4;7)$\vspace{-10pt}.}  
  \label{table:lower}
 \end{table}
Here $\Blue{a},\Blue{b},\Magenta{c},\Magenta{d}$ refer to the four conditions 
$(\Blue{\sigma_3}, \Blue{1362}, \, \Magenta{\sigma_4}, \Magenta{1346})$.
There are four other necklaces giving a monotone choice of points, and for those
the solutions were always real.
None of the remaining 8 necklaces had fewer than four real solutions.

Such lower bounds on the number of real solutions to enumerative geometric
problems were first found by Eremenko and Gabrielov~\cite{EG01b} in the context
of the Shapiro conjecture for Grassmannians.
Lower bounds have also been proven for problems of enumerating
rational curves on surfaces~\cite{IKS,Mi,W} 
and for some sparse polynomial systems~\cite{SS}.
We do not yet know a reason for the lower bounds here.

\subsubsection{Apparent upper bounds}
On $\Fl(1,2,3,4;5)$, set $\Blue{A}:=1325$ and $\Magenta{B}:=2143$.
The Schubert problem 
$(\Blue{A}^2, \Magenta{B}^3)$ has degree 7, but none of the 1 million
instances we computed had more than 5 real solutions.
 \begin{table}[htb]
  \begin{tabular}
    {|c||c|c|c|c|}\hline
     {Necklace} & \multicolumn{4}{c|}{Number of Real Solutions\rule{0pt}{13pt}}\\
     \cline{2-5}
        &1&3&5&7\rule{0pt}{13pt}\\\hline\hline
  $\Blue{AA}\Magenta{BBB}$ & 0 &500000 & 0 &   0   \\ \hline 
  $\Blue{A}\Magenta{B}\Blue{A}\Magenta{BB}$ & 193849 & 268969 & 37182 & 0 \\ \hline
  \end{tabular} \vspace{5pt}
  \caption{The Schubert problem $(\Blue{A}^2, \Magenta{B}^3)$ on $\Fl(1,2,3,4;5)$.}
  \label{T:upper}
 \end{table}

Neither condition $\Blue{A}$ nor $\Magenta{B}$  is Grassmannian, and so this
Schubert problem is not related to the conjectures in this paper.

\subsubsection{Apparent gaps}
On $\Fl(1,3,5;6)$, set $\Blue{A}:=21436$ and
$\Magenta{B}:=31526$.
The Schubert problem $(\Blue{A}^2, \Magenta{B}, \Red{\sigma_3}^2)$ 
has degree 8 and it appears to exhibit gaps in the possible numbers of real solutions.
Table~\ref{T:gaps} gives the data from this computation.
In each necklace, \Red{3}
represents the Grassmannian condition $\Red{\sigma_3}$. 
 \begin{table}[htb]
  \begin{tabular}
   {|c||c|c|c|c|c|}\hline
     {Necklace} & \multicolumn{5}{c|}{Number of Real Solutions\rule{0pt}{13pt}}\\
     \cline{2-6}
        &0&2&4&6&8\rule{0pt}{13pt}\\\hline\hline
 \Blue{AA}\Magenta{B}\Red{33}& 0  &     0& 991894 &0 &8106\\\hline  
 \Blue{AA}\Red{3}\Magenta{B}\Red{3}&111808 &  0& 888040& 0& 152 \\\hline
 \Blue{A}\Red{3}\Blue{A}\Red{3}\Magenta{B} &311285 & 0& 681416 &0& 7299 \\\hline 
 \Blue{A}\Red{33}\Blue{A}\Magenta{B}& 884186&  0 &115814 & 0& 0  \\\hline
%
%
  \end{tabular} \vspace{5pt}
  \caption{The Schubert problem 
       $(\Blue{A}^2, \Magenta{B}, \Red{\sigma_3}^2)$ on $\Fl(1,3,5;6)$.}
  \label{T:gaps}
 \end{table}
This is a new phenomena first observed in some sparse polynomial
systems~\cite[\S~7]{SS}.

\subsubsection{Small Galois group}\label{S:Galois}
One unusual problem that we looked at was on the flag manifold $\Fl(2,4;6)$
and it involved four identical non-Grassmannian conditions,
$1425$.
We can prove that this problem has six solutions, and that they are always all real.

\begin{thm}\label{T:sixReal}
  For any distinct $s,t,u,v\in\RP^1$, then intersection
\[ 
   X_{1425}(s) \cap X_{1425}(t) \cap X_{1425}(u) \cap X_{1425}(v)  
\]
 is transverse and consists of $6$ real points.
\end{thm}

This Schubert problem exhibits some other exceptional geometry concerning its Galois
group, which we now define.
Let $(w_1,\dotsc,w_s)$ be a Schubert problem for $\Flan$ and consider the
configuration space of $s$-tuples of flags $(F_\bullet^1,F_\bullet^2,\dotsc,F_\bullet^s)$
for which 
 \[
   X\ :=\ X_{w_1}F_\bullet^1 \cap X_{w_2}F_\bullet^2 \cap \dotsb\cap X_{w_s}F_\bullet^s
 \]
is transverse and hence $X$ consists of finitely many points.
If we pick a basepoint of this configuration space and follow the intersection
along a based loop in the configuration space, we will obtain a permutation of
the intersection $X$ corresponding to the base point.
Such permutations generate the \Blue{{\it Galois group}} of this Schubert
problem. 

Harris~\cite{Ha79} defined Galois groups for any enumerative geometric problem and
Vakil~\cite{Va04} investigated them for Schubert problems on Grassmannians,
showing that  many problems have a Galois group that contains at least the
alternating group. 
He also found some Schubert problems on Grassmannian whose Galois group is not
the full symmetric group. 
This Schubert problem also has a strikingly small Galois group, and is the
simplest Schubert problem we know with a small Galois group.

\begin{thm}\label{T:SmGal}
  The Galois group of the Schubert problem $(1425)^4$ on $\Fl(2,4;6)$ is the
  symmetric group on $3$ letters.
\end{thm}

We prove both theorems.
First, consider the Schubert variety $X_{1425}F_\bullet$
 \[
  X_{1425}F_\bullet\ =\ \{E_2\subset E_4\mid
   \dim E_2\cap F_3\geq 1\mbox{ and }\dim E_4\cap F_3 \geq 2\}.
 \]
The image of $X_{1425}F_\bullet$ under the projection $\pi_4\colon
\Fl(2,4;6)\twoheadrightarrow\Gr(4,6)$ is 
 \[
  \Omega_{1245}F_\bullet\ :=\ \{E_4 \in \Gr(4,6) \mid \dim E_4\cap F_3 \geq 2\}.
 \]
Since this Schubert variety has codimension 2 in $\Gr(4,6)$, a variety of dimension 8,
there are finitely many 4-planes $E_4$ which have Schubert position $1245$ with
respect to four general flags.
In fact, there are exactly 3.
(See Section~8.1 of~\cite{So97a}, which treats the dual problem in $\Gr(2,6)$.)

Thus we have a fibration of Schubert problems
 \begin{equation}\label{Eq:fibre}
   \bigcap_{i=1}^4 X_{1425}F_\bullet^i\ \xrightarrow{\ \pi_4\ }\ 
   \bigcap_{i=1}^4 \Omega_{1245}F_\bullet^i\,.
 \end{equation}
Let $K$ be a solution to the Schubert problem in $\Gr(4,6)$.
We ask, for which 2-planes $H$ in $\C^6$ is the flag $H\subset K$ a solution to the
Schubert problem in $\Fl(2,4;6)$?
From the description of $X_{1425}F_\bullet$, $H$ must be a 2-plane in $K$ which
meets each linear subspace $K\cap F^i_3$ non-trivially.
As $K$ lies in each Schubert cell $\Omega^\circ F_\bullet^i$, $K\cap F_3^i$ is a 2-plane.
Thus we are looking for the 2-planes $H$ in $K$ which meet four general 2-planes 
$K\cap F_3^i$.
There are two such 2-planes $H$, as this is an instance of the problem of lines in
$\P^3$ meeting four lines.
We conclude that there are six solutions to the Schubert problem on $\Fl(2,4;6)$.

This Schubert problem projects to one in $\Gr(2,6)$ with three solutions that
is dual to the projection in $\Gr(4,6)$.
Let $H_i$ and $K_i$ for $i=1,2,3$ be the 2-planes and 4-planes which are solutions to
the two projected problems.
For each $K_i$ there are exactly two $H_j$ for which
$H_j\subset K_i$ is a solution to the original problem in $\Fl(2,4;6)$.
Dually, for each $H_i$ there are exactly two $K_j$ for which 
$H_i\subset K_j$ is a solution to the original problem.
There is only one possibility for the configuration of the six flags, up to
relabeling: 
 \begin{equation}\label{E:crown}
  \raisebox{-20pt}{
   \begin{picture}(80,50)
   \put(4.5,4){\includegraphics[height=40pt]{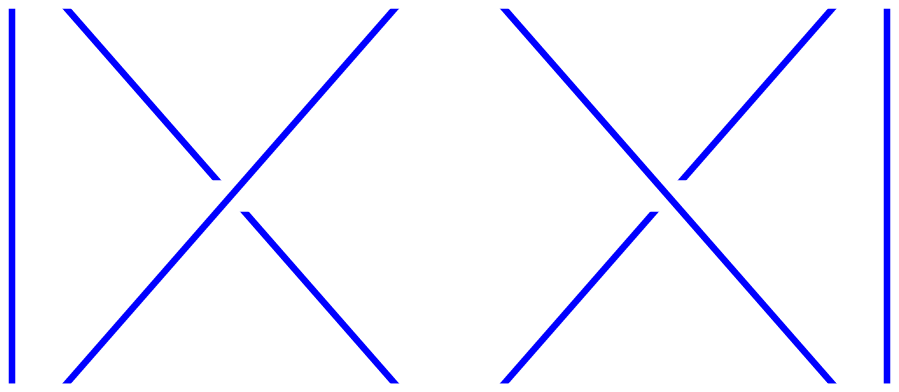}}
   \put(0,40){$K_3$}  \put(33,40){$K_2$} \put(66,40){$K_1$}
   \put(0, 0){$H_1$}  \put(33, 0){$H_2$} \put(66, 0){$H_3$}
  \end{picture}
}
 \end{equation}

\begin{proof}[Proof of Theorem~$\ref{T:sixReal}$]
 Since the flags osculate the rational normal curve, the
 problems obtained by projecting the intersection in Theorem~\ref{T:sixReal}
 to Grassmannians have only real solutions, as shown in Theorem~3.9
 of~\cite{So00a}. 
 Thus all subspaces $H_i$ and $K_i$ in~\eqref{E:crown}
 are real, and so the six solution flags of~\eqref{E:crown} are all real.
\end{proof}

\begin{proof}[Proof of Theorem~$\ref{T:SmGal}$]
 Since the six solution flags have the configuration given  in~\eqref{E:crown},
 we see that any permutation of the six solutions is 
 determined by its action on the three 4-planes $K_1, K_2, K_3$.
 Thus the Galois group is at most the symmetric group $S_3$.
 The explicit description given in Section~8.1 of~\cite{So97a} and also the
 analysis of Vakil~\cite{Va04} shows that the Galois group of the projected problem in
 $\Gr(4,6)$ is $S_3$.
\end{proof}

\subsubsection{A non-transverse Schubert problem}\label{S:improper}
Our experimentation uncovered a Schubert problem whose corresponding
intersection is not transverse or
even proper, when it involves flags osculating a rational normal curve.
This  may have negative repercussions for part of
Varchenko's program on the Bethe Ansatz and Fuchsian equations~\cite{MV04}. 
This was unexpected, as 
 Eisenbud and Harris showed that on a Grassmannian, any intersection
 \begin{equation}\label{Eq:genInt}
   X_{w_1}(t_1)\cap \dotsb\cap X_{w_m}(t_m)
 \end{equation}
 is proper in that it has the expected dimension
 $\dim(\alpha)-\sum\ell(w_i)$, if the points $t_1,\dotsc,t_m$ in $\P^1$ are
 distinct~\cite[Theorem 2.3]{EH83}. 
 On any flag manifold, if each condition (except possibly one) has codimension
 1 ($\ell(w_i)=1$), and if the points  $t_1,\dotsc,t_m\in\P^1$ are
 general, then the intersection~\eqref{Eq:genInt} is transverse, and hence
 proper~\cite[Theorem~2.1]{So99a}. 
 We show this is not the case for all Schubert problems on the flag manifold.

 The manifold of flags of type $\{1,3\}$ in $\C^5$ has dimension 8.
 Since $\ell(32514)=5$ and $\ell(21435)=2$, there are no flags of 
 type $\{1,3\}$ 
 satisfying the Schubert conditions $(325,\,(214)^2)$ imposed
 by three general flags.
 This is not the case if the flags osculate a rational normal curve $\gamma$.

\begin{thm}
 The intersection $X_{325}(u)\cap X_{214}(s)\cap X_{214}(t)$ is nonempty
 for all $s,t,u\in\P^1$.
\end{thm}

\begin{proof}
 We may assume without any loss that $u=\infty$, so that flags in
 $X^\circ_{325}(u)$ are given by matrices in $\calM_{325}$.
 Consider the $3\times 5$ matrix in  $\calM_{325}$.
 \begin{equation}\label{Eq:M32512}
   \left[ \begin{matrix}
     0&0&1&\frac{3}{2}(s+t)&6st\\
     0&1&0&-3st&0\\
     0&0&0&0&1\end{matrix}
   \right]
 \end{equation}
 Let $E_\bullet\colon E_1\subset E_3$ be the corresponding flag.
 We will show that $E_\bullet\in X_{214}(s)\cap X_{214}(t)$.
 Let $v_1$,  $v_2$, and $v_3$ to be the row vectors in~\eqref{Eq:M32512}.
 Consider the dual vector 
\[
   \lambda(s)\ :=\ (s^4,\,-4s^3,\,6s^2,\,-4s,\,1)\,,
\]
 and note that $\lambda(s)$ annihilates $\gamma(s)$, $\gamma'(s)$,
 $\gamma''(s)$, and $\gamma'''(s)$, so that $\lambda(s)$ is a linear form
 annihilating the 4-plane $F_4(s)$ osculating the rational normal curve
 $\gamma$ at the point $\gamma(s)$.
 Note that $v_1\cdot \lambda(s)^t=0$, so that $E_1\subset F_4(s)$.
 Also,
\[
   \gamma'(s)\ =\ 
    v_2 +2sv_1 +(4s^3-12s^2t)v_3\,,
\]
 and so $E_3\cap F_2(s)\neq 0$.
 In particular this implies that $E_\bullet\in X_{214}(s)$.
 We similarly have that $E_\bullet\in X_{214}(t)$.
\end{proof}

\section{Discussion}\label{S:Discussion}

We establish relationships between the different conjectures of Section 2,
between the conjectures for different Schubert problems on the same flag
manifold, and between the conjectures for Schubert problems on different flag
manifolds. 
This includes a proof of Theorem~\ref{T:simple} and a subtle generalization of
Conjecture~\ref{C:Main}. 
We conclude by proving Conjecture~\ref{C:Discriminant} for several Schubert
problems.

\subsection{Child problems}\label{S:T_simple}
The Bruhat order on $W^{\alpha}$ is defined by its covers $w\lessdot u$: if
$\ell(w)+1=\ell(u)$ and $w^{-1}u$ is a transposition $(b,c)$.
Necessarily, there exists an $i$ such that $b\leq \alpha_i<c$, but this number $i$
may not be unique. 
Write $w\lessdot_i u$ when $w\lessdot u$ in the Bruhat order and the
transposition $(b,c):=w^{-1}u$ satisfies $b\leq \alpha_i<c$.
This defines the cover relation in a partial order \Blue{$<_i$} on $W^{\alpha}$,
which is a subposet of the Bruhat order, and is called
the $\alpha_i$-Bruhat order in the combinatorics literature~\cite{So96}.
When $w<u$ are two Grasmannian permutations with the same descent $\alpha_i$ which
are related in Bruhat order, then $w<_iu$ and there is a chain of covers in the
$<_i$-order connecting $w$ to $u$.

Suppose that $(v,w_1,w_2,\dotsc,w_m)$ is a Schubert problem for $\Flan$
and that $v=\sigma_{\alpha_i}$.
For any permutation $u$ with $w_1\lessdot_i u$, we have
$\ell(v)+\ell(w_1)=\ell(u)$ and so 
$(u,\,w_2,\dotsc,w_m)$ is a Schubert problem for $\Flan$.
We say that $(u,\,w_2,\dotsc,w_m)$ is a \Blue{{\em child problem}} of the
original Schubert problem $(v,w_1,w_2,\dotsc,w_m)$ and write
\[
    (\Blue{v},\Blue{w_1},w_2,\dotsc,w_m)\ 
   \precdot\ (\Blue{u},\,w_2,\dotsc,w_m)\,,
\]
which defines the covering relation for a partial order $\prec$ on the set of
Schubert problems for $\Flan$.
Since every cover $w\lessdot u$ in the Bruhat order on $W^{\alpha}$ has the form
$\lessdot_i$ for some $i$, the minimal elements in this partial order $\prec$
are exactly the simple Schubert problems.
The reason for these definitions is the following lemma.

\begin{lemma}\label{L:induction}
 Suppose that 
 $(v,w_1,w_2,\dotsc,w_m)\precdot(u,\,w_2,\dotsc,w_m)$ is a cover between two
 Grassmannian Schubert problems for $\Flan$, where $\delta(w_1)=\alpha_i$,
 $v=\sigma_{\alpha_i}$, and $w_1\lessdot_i u$.
 If Conjecture~$\ref{C:Main-NT}$ holds for $(v,w_1,w_2,\dotsc,w_m)$, 
 then it holds for $(u,\,w_2,\dotsc,w_m)$.
\end{lemma}

The case $\beta=\alpha$ of Theorem~\ref{T:simple} follows from
Lemma~\ref{L:induction} as any Grassmannian Schubert problem is connected to a
simple Schubert problem via a chain of covers as in Lemma~\ref{L:induction}.
In turn, Lemma~\ref{L:induction} is a consequence of 
Lemma~\ref{L:limit}, which is proven in the next section.

\subsection{Limits of Schubert intersections}
Let $w\in W^{\alpha}$ be a Schubert condition for $\Flan$ and suppose that
$v=\sigma_{\alpha_i}$.
If $t\neq 0$, then the intersection $X_w(0)\cap X_v(t)$ is (generically)
transverse. 
One result of~\cite{So00b} concerns the limit of this intersection.
Specifically, we have the cycle-theoretic equality
 \begin{equation}\label{Eq:cycle}
   \lim_{t\to 0} X_w(0)\cap X_v(t)\ =\ 
    \sum_{w\lessdot_i u} X_u(0)\,.
 \end{equation}
That is, the support of the scheme-theoretic limit is the union of Schubert
varieties in the sum, and this scheme-theoretic limit is reduced at the generic
point of each Schubert variety in the sum.
We use this to prove the following lemma.

\begin{lemma}\label{L:limit}
 Let $(v,w_1,w_2,\dotsc,w_m)$ be a Schubert problem for $\Flan$, where 
 $v=\sigma_{\alpha_i}$.
 Suppose that $t_2,\dotsc,t_m$ are negative real numbers such that the
 intersection  
\[
    X_v(t)\cap X_{w_1}(0)\,\cap\, X_{w_2}(t_2)\cap\dotsb\cap X_{w_m}(t_m)
\]
 consists only of real points, for any positive number $t$.
 Then, for any permutation $u$ with $w_1\lessdot_i u$, the intersection
\[
    X_u(0)\,\cap X_{w_2}(t_2)\cap\dotsb\cap X_{w_m}(t_m)
\]
 consists only of real points
\end{lemma}

\begin{proof}
 Set $Y:=X_{w_2}(t_2)\cap\dotsb\cap X_{w_m}(t_m)$.
 We assumed that if $0<t$, then 
 $X_{w_1}(0)\cap X_v(t)\,\cap Y$ consists only of real points.
 The property of only having real points of intersection is preserved under
 taking limits, and so~\eqref{Eq:cycle} implies
 that every point of 
\[
   Y\ \cap\ \sum_{w_1\lessdot_i u} X_u(0)
\]
 is real.
 In particular, if $w_1\lessdot_i u$, then $Y\cap X_u(0)$ consists only of real
 points. 
\end{proof}

\begin{proof}[Proof of Lemma~$\ref{L:induction}$]
 Let $t_1,\dotsc,t_m\in\R\P^1$ be a monotone choice of points for 
 the Schubert problem $(u,w_2,\dotsc,w_m)$.
 Applying a real M\"obius transformation if necessary, we may assume that
 $t_1=0$ and that $t_2,\dotsc,t_m$ are negative real numbers.
 Thus it suffices to show that 
 \begin{equation}\label{Eq:limiting}
    X_u(0)\,\cap X_{w_2}(t_2)\cap\dotsb\cap X_{w_m}(t_m)
 \end{equation}
 consists only of real points.
 Since $\delta(u)=\delta(w_1)=\delta(v)=\alpha_i$, it follows that if $0<t$, then 
 $(t,0,t_2,\dotsc,t_m)$ is monotone with respect to the Schubert problem
 $(v,w_1,w_2,\dotsc,w_m)$.
 By our assumption that Conjecture~\ref{C:Main-NT} holds for
 $(v,w_1,w_2,\dotsc,w_m)$, the intersection 
\[
    X_v(t)\cap X_{w_1}(0)\,\cap\, X_{w_2}(t_2)\cap\dotsb\cap X_{w_m}(t_m)
\]
 consists only of real points, for any positive number $t$.
 But then Lemma~\ref{L:limit} implies that the intersection~\eqref{Eq:limiting}
 consists only of real points.
\end{proof}

\subsection{Refined monotone conjecture}
Lemma~\ref{L:limit} leads to an extension of Conjecture~\ref{C:Main}
to some cases when the Schubert problem is not Grassmannian.
We first give an example, which indicates a strengthening of Theorem~\ref{T:simple}.

\begin{ex}
 Consider the following instance of the cycle-theoretic
 equality~\eqref{Eq:cycle},
 \begin{equation}\label{Eq:14325}
  \lim_{x\to 0^+} X_{142}(0)\cap X_{\sigma_3}(x)\ 
    =\   X_{152}(0) \cup  X_{143}(0)\,. 
 \end{equation}
 Note that $\delta(142)=2$.
 Suppose that Conjecture~\ref{C:Main} holds for the Schubert problem
 $(\sigma_2^3,142,\sigma_3^3)$.
 Then, if $s<t<u<0<x<y<z$, the intersection
\[
   X_{\sigma_2}(s)\cap X_{\sigma_2}(t)\cap X_{\sigma_2}(u)\ \cap\
     X_{142}(0)\cap X_{\sigma_3}(x)\ 
    \cap   X_{\sigma_3}(y)\cap   X_{\sigma_3}(z)\,
\]
 consists only of real points, as the choice of points $s,t,u,0,x,y,z$ is
 monotone with respect to the given Schubert problem.
 As in the proof of Lemma~\ref{L:limit}, the limit~\eqref{Eq:14325} implies that
 whenever $s<t<u<0<y<z$, the intersection
\[
   X_{\sigma_2}(s)\cap X_{\sigma_2}(t)\cap X_{\sigma_2}(u)\ 
      \cap\ X_{143}(0)\ \cap X_{\sigma_3}(y)\cap   X_{\sigma_3}(z)\,
\]
 consists only of real points, even though the permutation $14325$ is not
 Grassmannian. 
\end{ex}

 We extend our notion of monotone choices of points to encompass this 
 last example.
 For a permutation $w\in W^{\alpha}$, let
 $\delta(w)\subset\{\alpha_1,\dotsc,\alpha_k\}$ be its 
 set of descents.
 Given two subsets $S, T\subset\{\alpha_1,\dotsc,\alpha_k\}$, we say
 that $S$ \Blue{{\it preceeds}} $T$, written $S<T$ if we have $i\leq j$ for all
 $i\in S$ and $j\in T$.
 This does not define a partial order on the set of subsets, but it does give a
 notion of when a list of subsets is increasing.
 For example
 \begin{equation}\label{Eq:new_monotone}
   \{2\}\ <\ \{2\}\ <\ \{2\}\ <\ \{2,3\}\ <\ \{3\}\ <\ \{3\}
 \end{equation}
 is increasing, but $\{2,3\}\not<\{2,3\}$.
Note that  $\{2\}<\{2\}$.

 A list of points $(t_1,\dotsc,t_m)\in\R\P^1$ is monotone with respect to
 a Schubert problem $(w_1,\dotsc,w_m)$ for $\Flan$ if the function 
\[
  t_i\ \longmapsto\ \delta(w_i)\subset \{\alpha_1,\dotsc,\alpha_k\}
\]
 is monotone, when the ordering of the $t_i$ is consistent with some ordering of
 $\R\P^1$. 
 For example, $(s<t<u<0<y<z)$ is monotone with respect to the Schubert problem 
 $(\sigma_2,\sigma_2,\sigma_2,143,\sigma_3,\sigma_3)$,
 as $\delta(143)=\{2,3\}$, and we have~\eqref{Eq:new_monotone}.
 We give a refinement of Conjecture~\ref{C:Main}, which drops the condition that
 the Schubert problem is Grassmannian.

\begin{conj}\label{C:New_Main}
 Suppose that $(w_1,\dotsc,w_m)$ is a Schubert problem for $\Flan$.
 Then the intersection 
 \begin{equation}\label{New_AA}
     X_{w_1}(t_1)\cap X_{w_2}(t_2)\cap 
    \dotsb\cap X_{w_m}(t_m)\,,
 \end{equation}
 is transverse with all points of intersection real, if the points
 $t_1,\dotsc,t_m\in\RP^1$ are monotone with respect to $(w_1,\dotsc,w_m)$.
\end{conj}

\begin{rem}
  There are many Schubert problems for which there are no monotone points.
  For example, two of the conditions $(\Blue{A})$ in the Schubert problem of
  Table~\ref{T:upper} have descent set $\{2,4\}$ and so there are no 
  monotone points.
  As reported there, for each of the two different necklaces, there are
  choices of points with not all  solutions real.
  Similarly, in the Schubert problem of Table~\ref{T:gaps}, there are three
  permutations with descent set $\{1,3,5\}$, and thus no monotone points.
  The Schubert problem in Section~\ref{S:Galois} consists of four
  identical conditions $w$ with $\delta(w)=\{2,4\}$, and so there are no
  monotone points.
  Nevertheless, we showed that all solutions are real.
\end{rem}

 The other conjectures of Section~\ref{S:conjectures} may be refined
 to include this more general notion of monotone points.
 For example, we conjecture that the discriminant of a Schubert problem does not
 vanish for monotone points, and that it (or its negative) lies in the
 preorder generated by differences of the $t_i$, as in 
 Conjecture~\ref{C:Discriminant}.

 The theorems of Section~\ref{S:conjectures} also hold in this generality, as the
 proofs are identical.
 For example, we have the following strengthening of
 Theorem~\ref{T:simple}.\medskip

\noindent{\bf Theorem~\ref{T:simple}$'$.}\  
{\it 
 Suppose that Conjecture~$\ref{C:New_Main}$ holds for all simple Schubert
 problems on a given flag manifold, $\Flan$.
 Then Conjecture~$\ref{C:New_Main}$ holds for all Schubert
 problems on any flag manifold 
 $\Fl(\beta,n)$ where $\beta$ is any subsequence of $\alpha$.
 (Here, the condition of transversality in Conjecture~$\ref{C:New_Main}$ is
 dropped.) 
}\medskip

\begin{ex} 
 Table~\ref{T:seven} shows data from the Schubert problem
 $(\Purple{\sigma_2}^2,\Blue{1432},
 \ForestGreen{135\,2},\Red{1254},\Magenta{\sigma_4}^2)$ on $\Fl(2,3,4;6)$, which has
 12 solutions, and involves two non-Grassmannian conditions.  
 In the necklaces, $\Purple{2},\Blue{A},\ForestGreen{3},\Red{B},\Magenta{4}$
 represent the five Schubert conditions, respectively.
 Their descent sets are 
 $\Purple{\{2\}},\Blue{\{2,3\}},\ForestGreen{\{3\}},\Red{\{3,4\}},\Magenta{\{4\}}$,
 \begin{table}[htb]
  \begin{tabular}
    {|c||c|c|c|c|c|c|c|}\hline
     {Necklace} & \multicolumn{7}{c|}{Number of Real Solutions\rule{0pt}{13pt}}\\
     \cline{2-8}
        &0&2&4&6&8&10&12\rule{0pt}{13pt}\\\hline\hline
  \Purple{22}\Blue{A}\ForestGreen{3}\Red{B}\Magenta{44} 
                   &0&0&  0&   0&    0&    0& 7500\\\hline
  \Purple{22}\Blue{A}\Red{B}\Magenta{44}\ForestGreen{3}
                   &0&0&  0&   0&    0&    0& 7500\\\hline
  \Purple{22}\Blue{A}\Red{B}\ForestGreen{3}\Magenta{44} 
                   &0&0&  0&   0&  306&  3776& 3416\\\hline
  \Purple{22}\Red{B}\ForestGreen{3}\Blue{A}\Magenta{44}
                  &0&0&  0&12& 1359& 3446&  2683  \\\hline
  \Purple{22}\ForestGreen{3}\Magenta{44}\Blue{A}\Red{B}
                  &0&0&  0&1213& 2129& 1771&  2387  \\\hline\hline
 
  \end{tabular} \vspace{5pt}
  \caption{The Schubert problem 
       $(\Purple{\sigma_2}^2,\Blue{1432},\ForestGreen{1354},\Red{1254},
          \Magenta{\sigma_4}^2)$ on $\Fl(2,3,4;6)$.}
  \label{T:seven}
 \end{table}
 so only the first row is monotone, and these data support
 Conjecture~\ref{C:New_Main}. 
 We only show 5 of the 90 necklaces.
\end{ex}

\subsection{Projections}\label{S:projections}
Suppose that  $\beta$ is a subsequence of $\alpha$.
In Section~\ref{S:basics} we considered projections 
$\pi\colon\Flan\to\Fl(\beta;n)$ obtained by forgetting 
the components of a flag $E_\bullet\in\Flan$ with dimension in 
$\alpha\setminus \beta$.
The image $\pi(X_wF_\bullet)$ of a Schubert variety of $\Flan$ 
is a Schubert variety of $\Fl(\beta;n)$ for a (possibly) different
permutation $\pi(w)$.
Recall that $w\in W^{\alpha}$ is a permutation whose descents can only occur at
positions in $\alpha$. 
The permutation $\pi(w)$ is obtained by ordering the values of $w$ between
successive positions in $\beta$.
For example, if $n=9$, $\alpha=\{2,4,5,7\}$ and $\beta=\{2,7\}$, then 
\[
  \pi(13\,\underline{58\,4\,27}\,69)\ =\ 13\,\underline{24578}\,69
   \qquad\mbox{and}\qquad
  \pi(26\,\underline{45\,7\,19}\,36)\ =\ 26\,\underline{14579}\,36\,.
\]

Because $\pi(X_wF_\bullet(s))=X_{\pi(w)}F_\bullet(s)$, if we have a Schubert problem 
$(w_1,\dotsc,w_m)$ on $\Flan$ and $m$ general flags, then 
$\pi$ is a map between the intersections
 \begin{equation}\label{Eq:project}
  \pi\ \colon\  X_{w_1}(t_1) \cap \dotsb\cap X_{w_m}(t_m)
  \ \longrightarrow\ 
   X_{\pi(w_1)}(t_1) \cap \dotsb\cap X_{\pi(w_m)}(t_m)\ .
 \end{equation}
Suppose that both $(w_1,\dotsc,w_m)$ and $(\pi(w_1),\dotsc,\pi(w_m))$ are 
Schubert problems.
Then the map $\pi$ of~\eqref{Eq:project} is a fibration with finite fibres. 
If the two problems have the same same degree, then
$\pi$ is an isomorphism.
In that case, we say that $(\pi(w_1),\dotsc,\pi(w_m))$ is a
\Blue{{\it projection}} of $(w_1,\dotsc,w_m)$ and that  $(w_1,\dotsc,w_m)$ is a
\Blue{{\it lift of}} $(\pi(w_1),\dotsc,\pi(w_m))$.

\begin{thm}\label{T:Projections}
 Suppose that the Schubert problem  $w:=(w_1,\dotsc,w_m)$ on $\Flan$ 
 is a lift of the Schubert problem $\pi(w)=(\pi(w_1),\dotsc,\pi(w_m))$ 
 on $\Fl(\beta;n)$.
 If Conjecture~$\ref{C:New_Main}$ holds for $\pi(w)$ then it holds for $w$.
\end{thm}

\begin{proof}
 Suppose that the permutations in $w$ are ordered so that 
\[
  \delta(w_1)\ <\ \delta(w_2)\ <\ \dotsb \ <\ \delta(w_m)
\]
 and let $t_1<\dotsb<t_m$ be real numbers.
 Then $\delta(\pi(w_1))<\dotsb <\delta(\pi(w_m))$
 and our assumption on $\pi(w)$ implies that the right-hand
 intersection in~\eqref{Eq:project} consists only of real points. 
 Since the map $\pi$ in~\eqref{Eq:project} is an isomorphism, we
 conclude that the left-hand intersection in~\eqref{Eq:project} consists only of
 real points.
\end{proof}

\begin{ex}
 Projection and lifts relate Schubert problems in many ways.
 The Grassmannian Schubert problem 
 $w:=(4\,1235, 15\,234, 135\,24, 1345\,2, 12456\,)$ on $\Fl(1,2,3,4,5;6)$ 
 has degree 4 and and it projects to the Schubert problem 
 $(\sigma_3, 125, 135, 134,\sigma_3)$ on the Grassmannian $G(3,6)$, which also
 has degree 4.
 One may compute a discriminant (as in~\cite[\S 3E]{So00a}) to show that the
 Shapiro conjecture holds for this Schubert problem.
 But then every Shapiro-type intersection for $w$ has all solutions real, and thus
 Conjecture~\ref{C:Main} holds for $w$.
 More interestingly, the projection of $w$ to $\Fl(2,4;6)$ also has only real
 solutions.
 This is the problem $(14\,23, 15\,23, 1325,  1345, 1245)$ of degree 4.
 Since the conditions have descents $(\{2\},\{2\},\{2,4\},\{4\},\{4\})$, there
 is a monotone choice of points, and so Conjecture~\ref{C:New_Main} holds
 for this last Schubert problem. 
\end{ex}


We now complete the proof of Theorem~\ref{T:simple}, showing that if
Conjecture~\ref{C:New_Main} holds for all simple Schubert problems on $\Flan$,
then Conjecture~\ref{C:New_Main} holds for all Schubert problems on
$\Fl(\beta,n)$, for any subsequence $\beta$ of $\alpha$.
Here, we drop the claim of transversality in Conjecture~\ref{C:New_Main}.
The proof will involve Schubert problems $w=(w_1,\dotsc,w_m)$ on $\Flan$ such
that $\pi(w)=(\pi(w_1),\dotsc,\pi(w_m))$ is a Schubert problem on
$\Fl(\beta;n)$, where $\pi\colon\Flan\to\Fl(\beta;n)$ is the projection map.
When this happens and the problem $w$ has non-zero degree, we say that the
Schubert problem $w$ is \Blue{{\it fibred over}} $\pi(w)$. 
Note that we do not require the two problems to have the same degree.
While it is not the case that $\pi(w)$ is a Schubert problem on $\Fl(\beta;n)$
whenever $w$ is a Schubert problem on $\Flan$, it turns out that for every
Schubert problem $v$ on $\Fl(\beta;n)$, there are many Schubert problems $w$ on
$\Flan$ which are fibred over $v$, and the degree of $w$ is always a positive
multiple of the degree of $v$.
The geometry behind this is discussed, for instance, in~\cite{PuSo}.

Indeed, the fibre $Y$ of the projection $\pi\colon\Flan\to\Fl(\beta;n)$
is a (product of) flag manifolds.
The map $\pi\colon X_wF_\bullet\to X_{\pi(w)}F_\bullet$ is almost a fibre bundle.
The fibre over a general point of $ X_{\pi(w)}F_\bullet$ is a Schubert variety in
$Y$ whose indexing permutation is $\pi(w)^{-1}w$.
Then if the flags are in general position, then $\pi$ restricts to a fibration
 \begin{equation}\label{Eq:fibration}
   \pi\ \colon\ 
   X_{w_1}F_\bullet^1 \cap \dotsb\cap X_{w_m}F_\bullet^m\ \longrightarrow\ 
   X_{\pi(w_1)}F_\bullet^1 \cap \dotsb\cap X_{\pi(w_m)}F_\bullet^m
 \end{equation}
with fibre the Schubert intersection in $Y$ given
by $(\pi(w_1)^{-1}w_1,\dotsc,\pi(w_m)^{-1}w_m)$.

When a problem $w$ is fibred over a problem $v$, there may be conditions $w_i$ of
$w$ such that $\pi(w_i)=\iota$, the identity permutation.
This condition $\iota$ is \Blue{{\em trivial\/}} because $X_\iota=\Fl(\beta;n)$.
Two problems $v$ and $v'$ are {\it equivalent} if they differ only in 
trivial conditions. 

The full statement of Theorem~\ref{T:simple} is a consequence
of the following result and the version when $\beta=\alpha$ already proven. 

\begin{thm}\label{T:fibred}
 Suppose that $\beta$ is a subsequence of $\alpha$ and that $v$ is a
 simple Schubert problem for $\Fl(\beta;n)$.
 Then there is a simple Schubert problem $w$ for $\Flan$ such that if
 Conjecture~$\ref{C:New_Main}$ holds for $w$, then it holds for $v$.
\end{thm}

\begin{proof}
 Suppose that $w=(w_1,\dotsc,w_m)$ is a simple Schubert problem on $\Flan$, 
 and each $w_i$ is a simple transposition of the form
 $\sigma_{\alpha_j}$, for some $j$.
 Then
\[
   \pi(\sigma_{\alpha_j})\ =\ \left\{\begin{array}{lcl}
    \sigma_{\alpha_j}&\ &\mbox{\rm if $\alpha_j\in\beta$}\\
     \iota&&\mbox{\rm otherwise\,.}
  \end{array}\right.
\]
 It follows that $\pi(w)$ is a simple Schubert problem on $\Fl(\beta;n)$ which
 involves some trivial Schubert varieties $X_\iota$.
 As in the proof of Theorem~\ref{T:Projections}, if $(t_1,\dotsc,t_m)$ is
 monotone for $w$, then it will be monotone for $\pi(w)$.
 Note that if $\pi(w_i)=\iota$, then the choice of the point $t_i$ does not
 affect the Schubert intersection for $\pi(w)$.

 The converse is also true.
 Let $v$ be the Schubert problem $\pi(w)$, where we have dropped all of the
 trivial conditions $\iota$.
 Any monotone choice of points for $v$ may be extended to a
 monotone choice of points $(t_1,\dotsc,t_m)$ for $w$.
 We need only choose points $t_i$ for those $w_i$ such that $\pi(w)=\iota$
 in a way to preserve monotonicity, which is easy.

 Suppose now that $v$ is a simple Schubert problem on $\Fl(\beta;n)$.
 Then there is a simple Schubert problem $w$ on $\Flan$ which is fibred over
 $v$.
 Indeed,  let $Y$ be the flag manifold which is the fibre of the projection
 $\pi\colon\Flan\to\Fl(\beta;n)$.  
 It suffices to add simple Schubert conditions to $v$ coming from any
 simple Schubert problem on $Y$ with degree $>0$.
 These added conditions $w_i$ have descents in $\alpha\setminus\beta$, 
 so the Schubert problems $\pi(w)$ and $v$ are equivalent. 
 Pick a monotone choice of points for $v$ and, as explained in the previous
 paragraph, extend it to a monotone choice of points for $w$.
 If Conjecture~\ref{C:New_Main} holds for $w$, then all the points in 
 $X_{w_1}(t_1)\cap\dotsb\cap X_{w_m}(t_m)$ are real.
 The map $\pi$~\eqref{Eq:fibration} exhibits this as a surjection onto 
 $X_{\pi(w_1)}(t_1)\cap\dotsb\cap X_{\pi(w_m)}(t_m)$, which
 equals the corresponding intersection for the Schubert problem $v$ and the
 original monotone choice of points.
\end{proof}

\begin{ex}
 Theorem~\ref{T:fibred} involved one Schubert problem fibred over another.
 An example is provided by the Schubert problem
 $(\Blue{\sigma_2}^4, (\Magenta{1245})^4)$ on $\Fl(2,4;6)$, which has degree 6.
 As with the example in Section~\ref{S:Galois}, this is fibred over the Schubert 
 problem on $\Gr(4,6)$ involving the intersection of four Schubert varieties 
 $\Magenta{\Omega_{1245}}$ given by flags osculating the rational normal curve
 at the  points corresponding to the conditions $\Magenta{1245}$.
 All three solution 4-planes $K_1,K_2$, and $K_3$ are real, and the fibre over
 $K_i$ is the the problem of four 2-planes in $K_i$ meeting four two planes
 that are the intersection of $K_i$ with four 4-planes osculating the rational
 normal curve at the points corresponding to the conditions $\Blue{\sigma_2}$.

 This problem in the fibre is not equivalent to an instance of the Shapiro
 conjecture for 2-planes in the 4-space $K_i$.
 If it were equivalent to an instance of the
 Shapiro conjecture, then all solutions for every necklace of
 Table~\ref{T:fibration} would be real, which is not the case.
 In the necklaces of Table~\ref{T:fibration}, $\Blue{2}$ represents the condition
 $\Blue{\sigma_2}$ and $\Magenta{4}$ represents the condition
 $\Magenta{1245}$.
 \begin{table}[htb]
  \begin{tabular}
    {|c||c|c|c|c|}\hline
     {Necklace} & \multicolumn{4}{c|}{Number of Real Solutions\rule{0pt}{13pt}}\\
     \cline{2-5}
        &0&2&4&6\rule{0pt}{13pt}\\\hline\hline
   \Blue{2222}\Magenta{4444} 
    &  0 &    0 &     0 &100000\\\hline
   \Blue{222}\Magenta{4}\Blue{2}\Magenta{444}
    &  0 &    0 &     0 & 100000\\\hline
   \Blue{222}\Magenta{44}\Blue{2}\Magenta{44}
    &  0 &    0 &     0 & 100000\\\hline
   \Blue{22}\Magenta{44}\Blue{22}\Magenta{44}
    &  0 &    0 &   122 & 99878  \\\hline
  \Blue{22}\Magenta{4}\Blue{2}\Magenta{4}\Blue{2}\Magenta{44}
    &  0 &   12 &  3551 & 96437  \\\hline
  \Blue{22}\Magenta{4}\Blue{2}\Magenta{44}\Blue{2}\Magenta{4}
    &  0 &  105 &  8448 & 91447 \\\hline
  \Blue{2}\Magenta{4}\Blue{2}\Magenta{4}\Blue{2}\Magenta{4}\Blue{2}\Magenta{4}
    &  0 & 1050 & 19964 & 78986  \\\hline
  \Blue{22}\Magenta{4}\Blue{22}\Magenta{444}   
    & 18 &   340 &  5147 & 94495  \\\hline
  \end{tabular} \vspace{5pt}
  \caption{The Schubert problem 
       $(\Blue{\sigma_2}^4, (\Magenta{1245})^4)$ on $\Fl(2,4;6)$.}
  \label{T:fibration}
 \end{table}
%
%
%

\end{ex}

\subsection{Discriminants}\label{S:discriminants}
Let $w=(w_1.\dotsc,w_m)$ be a Schubert problem on $\Flan$.
The {\it discriminant} 
$\Sigma\subset (\P^1)^m$ is the set of points $(t_1,\dotsc,t_m)$ where the
intersection
 \begin{equation}\label{Eq:disc_int}
  X_{w_1}(t_1)\cap \dotsb \cap X_{w_m}(t_m)
 \end{equation}
is not transverse.
When $\Sigma\neq(\P^1)^m$, this is a hypersurface defined by the
discriminant polynomial $\Delta_w(t_1,\dotsc,t_m)$, which is separately
homogeneous in each homogeneous parameter $t_i$.
For each of three Schubert problems, we will prove a weaker version of
Conjecture~\ref{C:Discriminant} which implies Conjecture~\ref{C:Main}. 

This version is weaker because we do not compute
$\Delta_w(t_1,\dotsc,t_m)$, as that would be infeasible.
Instead, we will fix three  parameters, say $t_1=\infty$, $t_2=0$, and $t_3=1$
(or $t_3=-1$). 
Then we can carry out the computation in the local coordinates $\calM_{w_1}$ for
$X^\circ_{w_1}(\infty)$, or in local coordinates for the intersection of two
cells $X^\circ_{w_1}(\infty)\cap X^\circ_{w_2}(0)$.
In these coordiantes, we generate the ideal defining the
intersection~\eqref{Eq:disc_int}, compute an
eliminant $F(x;t)$ for one of our coordinates, and then
compute its discriminant $\Delta_x(t)$. 

If we specialize the parameters $t_1, t_2$, and $t_3$ to these fixed values,
then $\Delta_w(t_1,\dotsc,t_m)$ will divide $\Delta_x(t)$, but there may be
other factors in $\Delta_x(t)$.
We minimized these extraneous factors by computing the greatest common divisor
of these discriminants $\Delta_x(t)$ for each coordinate $x$.
We also remove factors common to a leading term of any eliminant,
as those correspond to solutions that are not on our chosen coordinate patch.

\begin{thm}\label{T:discriminants}
 Conjecture~$\ref{C:Discriminant}$ holds for the following Schubert problems.
\begin{enumerate}

 \item $(\Purple{\sigma_1}, \Blue{4\,3\,1256}, \ForestGreen{13\,25\,46},
        \Red{1256\,4\,3}, \Magenta{\sigma_5})$ on $\Fl(1,2,4,5;6)$.  
        This has $2$ solutions. 

 \item $(\Blue{24\,135}, \Blue{13\,245}, 
        \Magenta{134\,25}, (\Magenta{124\,35})^2)$ on $\Fl(2,3;5)$.   
        This has $3$ solutions. 

 \item $(\Blue{146\,2357}, \Blue{135\,2467},  
       \Magenta{1246\,357}, \Magenta{1256\,347})$ on $\Fl(3,4;7)$.   
        This has $4$ solutions.

\end{enumerate}

\end{thm}

\begin{proof} 
(1)\ 
 Consider the Schubert intersection
\[
   X_{\Purple{\sigma_1}}(t)\,\cap\,  
   X_{\Blue{4\,3\,1256}}(\infty)\,\cap\, 
   X_{\ForestGreen{13\,25\,46}}(-1)\,\cap\,  
   X_{\Red{1256\,4\,3}}(0)\,\cap\,  
   X_{\Magenta{\sigma_5}}(s)\,,
\]
 on $\Fl(1,2,4,5;6)$. 
 Since these Schubert conditions have respective descent sets
\[
   \Purple{\{1\}},\ \Blue{\{1,2\}},\ \ForestGreen{\{2,4\}},\ 
    \Red{\{4,5\}},\ \Magenta{\{5\}}\,,
\]
 the set of monotone points is $\{(s,t)\mid 0<s<t\}$.
 The discriminant we computed had two factors.
 One was $s^6$ and here is the other factor
 \begin{align*}
  &\quad\; 2500s^4t^4\ +\ 18000s^3t^4+4000s^4t^3\ +\ 
           50000s^2t^4+31100s^3t^3+2260s^4t^2\\
  & +64000st^4+91400s^2t^3+20040s^3t^2+480s^4t\\
  & +32000t^4+122800st^3+63905s^2t^2+5550s^3t+9s^4\\
  & +64000t^3+91400st^2+20040s^2t+480s^3\\
  & +50000t^2+31100st+2260s^2\ +\ 18000t+4000s\ +\ 2500\,.
 \end{align*}
 This is a positive sum of monomials and is thus positive
 when $0<s,t$, which includes the set of monotone points.
\smallskip

\noindent (2)\ 
 Consider the Schubert intersection 
\[
  X_{\Blue{24\,135}}(\infty)\,\cap\,
  X_{\Blue{13\,245}}(-1)\,\cap\,
  X_{\Magenta{134\,25}}(0)\,\cap\,
  X_{\Magenta{124\,35}}(s)\,\cap\,
  X_{\Magenta{124\,35}}(t)
\] 
 on the flag variety $\Fl(2,3;5)$.
 Since these Schubert conditions have respective descents at 
 $\Blue{2},\Blue{2},\Magenta{3},\Magenta{3},\Magenta{3}$, the
 set of monotone points is 
 \begin{equation}\label{Eq:st_monotone}
   \{(s,t)\mid -1<s,t,\ s,t\neq 0,\ \mbox{and } s\neq t\}\,.
 \end{equation}
 We display the discriminant, shading the region with monotone points.
\[
  \begin{picture}(170,170)
   \put(0,0){\includegraphics[height=170pt]{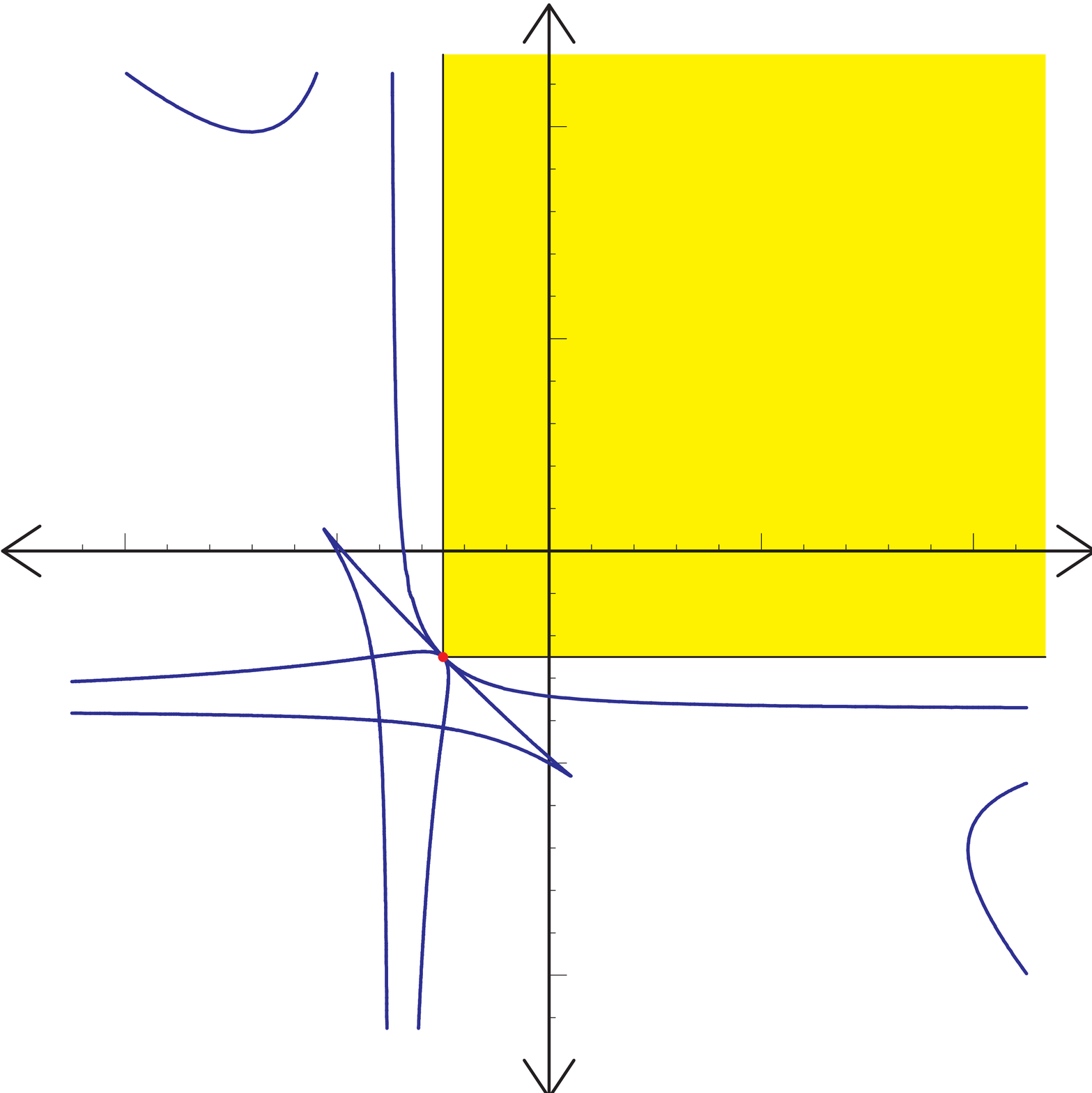}}
   \put(159,75){$s$}   \put(76,157){$t$}
   \put(  4,91){$-4$}  \put(91, 12.5){$-4$}
   \put( 39,91){$-2$}  \put(91, 47){$-2$}
   \put(117,91){$2$}   \put(92,116){$2$}
   \put(151,91){$4$}   \put(92,151){$4$}
  \end{picture}
\]

 In Figure~\ref{F:disc}, we write this discriminant in terms of $t$ and $s-t$, 
 whose positivity defines the region where $0<s<t$, a subset of the set of monotone
 points. 
 \begin{figure}
 \begin{align*}
  &\quad\;  800t^9+3600t^8(s{-}t)+7744t^7(s{-}t)^2+10304t^6(s{-}t)^3
           +8736t^5(s{-}t)^4+4480t^4(s{-}t)^5\\
  & +1248t^3(s{-}t)^6+144t^2(s{-}t)^7+5760t^8+23040t^7(s{-}t)
           +45792t^6(s{-}t)^2+56736t^5(s{-}t)^3\\
  & +43632t^4(s{-}t)^4+19584t^3(s{-}t)^5+4608t^2(s{-}t)^6
            +432t(s{-}t)^7+17888t^7+62608t^6(s{-}t)\\
  & +114816t^5(s{-}t)^2+130520t^4(s{-}t)^3+87520t^3(s{-}t)^4
                 +32064t^2(s{-}t)^5+5616t(s{-}t)^6\\
  & +324(s{-}t)^7+31712t^6+95136t^5(s{-}t)+161496t^4(s{-}t)^2+164432t^3(s{-}t)^3\\
  & +90048t^2(s{-}t)^4+23688t(s{-}t)^5+2268(s{-}t)^6+35456t^5+88640t^4(s{-}t)\\
  & +141256t^3(s{-}t)^2+123244t^2(s{-}t)^3+48726t(s{-}t)^4+6777(s{-}t)^5+25376t^4\\
  & +50752t^3(s{-}t)+79184t^2(s{-}t)^2+53808t(s{-}t)^3+11394(s{-}t)^4+10752t^3\\
  & +16128t^2(s{-}t)+27264t(s{-}t)^2+10944(s{-}t)^3+2048t^2
      +2048t(s{-}t)+4608(s{-}t)^2\ .
 \end{align*}
 \caption{A discriminant.\label{F:disc}}
 \end{figure}
 This discriminant is a positive linear combination of 49
 homogeneous monomials of degree 9 in the terms $t$ and $s-t$, and is thus
 positive on the region defined by $0<s<t$.
 There is a similar positive expression for the discriminant in
 terms of $1+t$ and $s$, and another in terms of $1+t$, $s-t$, and $-s$.
 Together with the expression in Figure~\ref{F:disc}, these show that 
 the discriminant is positive on the set $-1<s<t$ with  $s,t\neq 0$.
 Since the discriminant is symmetric in $s$ and $t$, the symmetric counterpart 
 of these three expressions shows that the discriminant is positive on the
 set~\eqref{Eq:st_monotone} of monotone points. 
\smallskip

\noindent(3)\ 
 Consider the Schubert intersection 
\[
 X_{\Blue{146\,2357}}(\infty)\,\cap\,
 X_{\Blue{135\,2467}}(0)\,\cap\,
 X_{\Magenta{1246\,357}}(s)\,\cap\,
 X_{\Magenta{1256\,347}}(1)
\]
on the flag variety $\Fl(3,4;7)$.
Since the Schubert conditions have descents
$\Blue{3},\Blue{3},\Magenta{4},\Magenta{4}$, the pointss are monotone when $0<s$.
Removing factors of $s$ and $1+s$ from the discriminant, we obtain 
\begin{align*}
  &\quad\;
   3515625+45243750s+221792500s^2+565872594s^3+777678231s^4\\
  &+1273923370s^5+932192307s^6+909742337s^{10}+1560886138s^{11}\\
  &+867109112s^{12}+367416324s^{13}+114976512s^{14}+13608000s^{15}+648000s^{16}\\
  &+\bigl(42966406s^3 + 352158344s^4 + 135425340s^5\bigr)(1-s^4)^2 \ ,
 \end{align*}
which is obviously positive when $0<s$.
\end{proof}

\begin{rem}
The first discriminant we computed, for the Schubert intersection
 \begin{equation}\label{Eq:F12456} 
   X_{\Purple{\sigma_1}}(\Purple{t})\,\cap\,  
   X_{\Blue{4\,3\,1256}}(\infty)\,\cap\, 
   X_{\ForestGreen{13\,25\,46}}(-1)\,\cap\,  
   X_{\Red{1256\,4\,3}}(0)\,\cap\,  
   X_{\Magenta{\sigma_5}}(\Magenta{s})\,,
 \end{equation}
was positive on more than just the monotone region.

Figure~\ref{Fig:table_discr} compares the table for this Schubert problem with a  
plot of the discriminant, which proves that 9 of the 12 necklaces will give only
real solutions.
\begin{figure}[htb]
\[
 \begin{picture}(420,219)
 \put(-5,107){
  \begin{tabular}
   {|c||c||c|c|}\hline
    Label& {Necklace}  &0&2\rule{0pt}{13pt}\\\hline\hline
  I   &\Blue{A}\ForestGreen{B}\Red{C}\Magenta{s}\Purple{t}&     0 & 100000\\\hline 
  II  &\Blue{A}\ForestGreen{B}\Red{C}\Purple{t}\Magenta{s}&     0 & 100000\\\hline 
  III &\Blue{A}\ForestGreen{B}\Magenta{s}\Purple{t}\Red{C}&     0 & 100000\\\hline 
  IV  &\Blue{A}\ForestGreen{B}\Purple{t}\Magenta{s}\Red{C}&     0 & 100000\\\hline 
  V   &\Blue{A}\Purple{t}\Magenta{s}\ForestGreen{B}\Red{C}&     0 & 100000\\\hline 
  VI  &\Blue{A}\Magenta{s}\Purple{t}\ForestGreen{B}\Red{C}&     0 & 100000\\\hline 
  VII &\Blue{A}\ForestGreen{B}\Magenta{s}\Red{C}\Purple{t}&     0 & 100000\\\hline 
  VIII&\Blue{A}\Purple{t}\ForestGreen{B}\Red{C}\Magenta{s}&     0 & 100000\\\hline 
  IX  &\Blue{A}\Purple{t}\ForestGreen{B}\Magenta{s}\Red{C}&     0 & 100000\\\hline 
  X   &\Blue{A}\ForestGreen{B}\Purple{t}\Red{C}\Magenta{s}& 24976 &  75024\\\hline 
  XI  &\Blue{A}\Magenta{s}\ForestGreen{B}\Red{C}\Purple{t}& 26065 &  73935\\\hline 
  XII &\Blue{A}\Magenta{s}\ForestGreen{B}\Purple{t}\Red{C}& 38023 &  61977\\\hline
  \end{tabular}\vspace{5pt}
}
  \put(195,-5){
   \begin{picture}(215,215)  
    \put(0,5){\includegraphics[height=220pt]{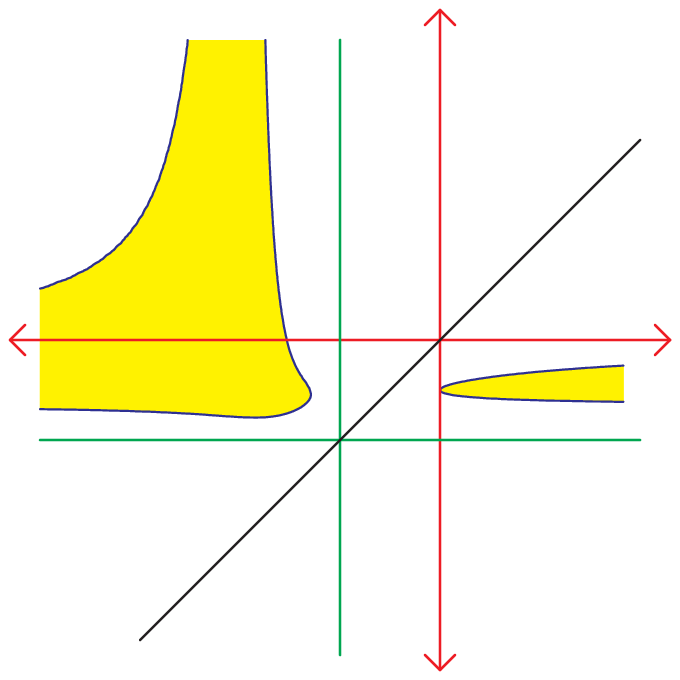}}
      \put(200,150){I}  \put( 87, 40){V}   \put(120, 40){IX}
      \put(180,165){II} \put( 67, 55){VI}  \put(216, 94){X}
      \put(129, 85){III}\put(120,150){VII} \put( 55,150){XI}
      \put(117,101){IV} \put(165, 40){VIII}\put( 53, 97){XII}
      \put( 76, 10){$\Purple{s}=\ForestGreen{B}$}
      \put(151, 10){$\Purple{s}=\Red{C}$}
      \put(180, 70){$\Magenta{t}=\ForestGreen{B}$}
      \put(180,120){$\Magenta{t}=\Red{C}$}
   \end{picture}}

  \end{picture}
\]
\caption{The Schubert problem $(\Purple{\sigma_1}, \Blue{4\,3\,1256}, 
   \ForestGreen{13\,25\,46}, \Red{1256\,4\,3},\Magenta{\sigma_5})$ on
   $\Fl(1,2,4,5;6)$ and its discriminant.}\label{Fig:table_discr}
\end{figure}
Indeed, the shaded region is where the discriminant is negative.
The $(\Magenta{s},\Purple{t})$-plane is divided into 12 regions by the lines
$\Magenta{s}=\Purple{t}$ and $\Magenta{s},\Purple{t}= \Red{0},\ForestGreen{-1}$, which
are points that cannot be used in the intersection~\eqref{Eq:F12456}. 
Each region corresponds to a necklace, and is labeled by the row of its
corresponding necklace.  
For the necklaces, we use $\Purple{t}$, $\Blue{A}$, $\ForestGreen{B}$, $\Red{C}$, and
$\Magenta{s}$  to denote the conditions
$\Purple{\sigma_1}$, $\Blue{4\,3\,1256}$, $\ForestGreen{13\,25\,46}$,
$\Red{1256\,4\,3}$, and $\Magenta{\sigma_5}$, respectively. 
\end{rem}

\section{Methods}\label{S:Methods}

The {\it raison d'\^{e}tre} for this paper is our computer experimentation
investigating the number of real points in Schubert intersections of the
form
 \begin{equation}\label{Eq:Schub_int_prob}
   X_{w_1}(t_1)\cap X_{w_2}(t_2)\cap \dotsb\cap X_{w_m}(t_m)\,,
 \end{equation}
for Schubert problems $(w_1,w_2,\dotsc,w_m)$ on small flag manifolds.
We determined this number for 520,420,135 different
intersections involving 1126 different Schubert problems on 29
different flag manifolds.
This used 15.76 gigahertz-years of computer time.

Table~\ref{table:conjs} shows the effort devoted to
studying the three main conjectures:
The Shapiro Conjecture for Grassmannians (Conjecture~\ref{C:SC}), our
Monotone Conjecture for Grassmannian Schubert problems
(Conjecture~\ref{C:Main}), and the Refined Monotone Conjecture
(Conjecture~\ref{C:New_Main}). 
Since these are in increasing order of generality, each of the last two rows of
Table~\ref{table:conjs} only shows the extra effort devoted to the
corresponding conjecture. 
 \begin{table}[htb]
  \begin{tabular}
   {|l||c|r|c|}\hline
   & Number of& \multicolumn{1}{c|}{Number of}& gigahertz\\
   & Problems& \multicolumn{1}{c|}{Intersections}&years\\\hline
  Conjecture~\ref{C:SC}      & 212&132,919,238& 3.57\\\hline
  Conjecture~\ref{C:Main}    & 376& 25,524,191& 1.23\\\hline
  Conjecture~\ref{C:New_Main}& 201&  7,223,660& 0.77\\\hline
 \end{tabular}\vspace{5pt}
  \caption{Resources devoted to the conjectures.} 
  \label{table:conjs}
 \end{table}
%
The numbers for the last two conjectures are only a small
fraction of the total effort expended in this experimentation.
This is because only a small fraction of necklaces are 
monotone.

A significant part of our investigation was devoted to
the Shapiro Conjecture for Grassmannians
(Conjecture~\ref{C:SC}), for $\Gr(3,6)$, $\Gr(3,7)$, and $\Gr(4,8)$.
While this conjecture had been studied before~\cite{So00a}, 
the scope of previous experiments was limited. 

Section 4.1 explains how we determined the number of real solutions in an 
intersection~\eqref{Eq:Schub_int_prob}.
Section 4.2 describes how we 
investigated such intersections for many necklaces and choices of
points for a single Schubert problem.
Section 4.3 discusses the design of the experiment, that is, how we chose which
Schubert problems to investigate. 

\subsection{Computation of a single Schubert intersection~(4.1)}
All computations were done on Intel processors running
Linux, using the computer algebra systems Singular {\tt 2-0-5}~\cite{SINGULAR} and
Maple, which were called from {\tt bash} shell scripts.
Maple managed the data, created the Singular scripts,
and counted the real solutions to univariate eliminants.

To study a Schubert intersection~\eqref{Eq:Schub_int_prob}, we
generated the ideal of the intersection in local coordinates
$\calM_{w_1}$ parametrizing the Schubert cell
$X^\circ_{w_1}(\infty)$. 
For this, we fixed $t_1=\infty$, and the other points
$t_2,\dotsc,t_m$ were rational numbers, 
and the ideal was generated by the equations for each Schubert
variety $X_{w_i}(t_i)$ as described in Section~\ref{S:basics} and
in Section~\ref{S:SC} (where the flags $F_\bullet(t_i)$ were described).
Because Gr\"obner basis computation is extremely sensitive to the
number of variables, the first Schubert condition $w_1$ was chosen to
minimize the number of coordinates in the parametrization 
$\calM_{w_1}$ of the Schubert cell $X^\circ_{w_1}(\infty)$.

Singular computed a degree reverse lexicographic Gr\"obner basis for
this ideal and then used the FGLM algorithm~\cite{FGLM} 
to compute a square-free univariate eliminant with degree equal to the
degree of the Schubert problem.
This guaranteed that the original intersection was transverse and that
its number of real points is equal to the number of real roots of the
eliminant (see the discussion in~\cite[\S 2.2]{So_M2}).
This number of real roots was computed using the {\tt realroot}
command in Maple.
When such an eliminant could not be computed, data describing the
intersection were set aside and later studied by hand.

\subsection{Investigation of a single Schubert problem}
For a given Schubert problem $(w_1,\dotsc,w_m)$, we determined the
number of real points in many different Schubert intersections of the
form~\eqref{Eq:Schub_int_prob}. 
Once a problem was selected, data necessary for the
experimentation were precomputed and stored in a data file.
These data included a list $L$ of permutations of the numbers $\{2,\dotsc,m\}$
and a set $S$ of rational numbers.
The list $L$ typically consisted of one permutation representing each necklace
we decided to investigate.
This data file was updated throughout the computation as it also recorded the 
numbers of real solutions found for the different necklaces and for
different choices of points.

Most  Schubert problems were run on a single computer.
The actual computation was organized by a shell script, whose main part was a
loop. 
In each iteration, the loop variable was used as a seed for Maple's 
random-number generator to select a random subset $t_2,\dotsc,t_m$ of the points 
from $S$, which were ordered so that $t_2<\dotsb<t_m$.
For each permutation $\sigma$ of $L$, the number of real points in the
intersection 
 \begin{equation}\label{E:permutation}
   X^\circ_{w_1}(\infty)\,\cap\, 
   X_{w_2}(t_{\sigma(2)})\cap X_{w_3}(t_{\sigma(3)})\cap 
   \dotsb\cap  X_{w_m}(t_{\sigma(m)})
 \end{equation}
was determined and included in the data file.
The data file also kept track
of the CPU time used in the computation, and recorded the average size of the
univariate eliminants.
The number of iterations of the shell script depended upon our
interest in the problem and the computational cost.

\begin{figure}
\begin{center}
\fbox{
\begin{minipage}{430pt}
\begin{flushleft}
{\large\bf Enumerative problem $\Purple{W_{\includegraphics{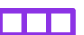}}}
 \Blue{(X_{\includegraphics{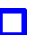}})^2}
 \ForestGreen{(Y_{\includegraphics{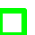}})^4}=7$  on $\Fl(1,2,3;5)$\rule{0pt}{17pt}}

\Blue{\rule{430pt}{0.5pt}}

\ \includegraphics[height=15pt]{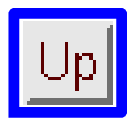}\rule{0pt}{17pt}\vspace{-10pt}

\Blue{\rule{430pt}{0.5pt}}

\begin{tabular}{ll}
 \multicolumn{1}{c}{\Blue{{\bf Experimental data}}}
 &
 \multicolumn{1}{c}{\Blue{{\bf Related Problems}}}\\
 \multirow{5}{215pt}{
  \begin{tabular}{||c||r|r|r|r||}\hline
   & \multicolumn{4}{c||}{Number of Real Solutions\rule{0pt}{13pt}}\\\hline
   {\small Necklace} &\multicolumn{1}{c|}{{\bf 1}}&\multicolumn{1}{c|}{{\bf 3}}
              &\multicolumn{1}{c|}{{\bf 5}}&\multicolumn{1}{c||}{{\bf 7}}\\\hline
   {\small{\bf\Purple{a}\Blue{bb}\ForestGreen{cccc}}}$\!$
          & 0&0&0&25000\\\hline
   {\small{\bf\Purple{a}\Blue{b}\ForestGreen{cccc}\Blue{b}}}$\!$ 
          &0&0&0&25000\\\hline
   {\small{\bf\Purple{a}\ForestGreen{cc}\Blue{bb}\ForestGreen{bb}}}$\!$ 
          &0&89&10500&14411\\\hline
   {\small{\bf\Purple{a}\ForestGreen{c}\Blue{bb}\ForestGreen{ccc}}}$\!$ 
          &0&2374&5740&16886\\\hline
   {\small{\bf\Purple{a}\Blue{b}\ForestGreen{cc}\Blue{b}\ForestGreen{cc}}}$\!$ 
          &0&2560&13204&9236\\\hline
   {\small{\bf\Purple{a}\Blue{b}\ForestGreen{ccc}\Blue{b}\ForestGreen{c}}}$\!$ 
          &0&4456&9753&10791\\\hline
   {\small{\bf\Purple{a}\ForestGreen{a}\Blue{b}\ForestGreen{cc}\Blue{b}\ForestGreen{c}}}$\!$ 
          &29&2571&14627&7773\\\hline
   {\small{\bf\Purple{a}\Blue{b}\ForestGreen{c}\Blue{b}\ForestGreen{ccc}}}$\!$ 
          &1120&5364&9633&8883\\\hline
   {\small{\bf\Purple{a}\ForestGreen{c}\Blue{b}\ForestGreen{c}\Blue{b}\ForestGreen{cc}}}$\!$ 
          &3446&5566&9132&6856\\\hline
  \end{tabular}}
  &\Red{Projections}\rule{0pt}{20pt}\\
  &\begin{tabular}{||l|l|l||}\hline
   Variety& Problem &$\!$\#$\!$\\\hline
   {\small \Blue{\underline{Fl(2,3;5)}}}\raisebox{-7pt}{\rule{0pt}{19pt}}&
   {\small $\Blue{(X_{\includegraphics[height=5pt]{figures/X1.eps}})^2
           X_{\includegraphics[height=5pt]{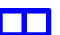}}}
     \ForestGreen{(Y_{\includegraphics[height=5pt]{figures/Y1.eps}})^4}$}& 7\\\hline
  \end{tabular}
  \\&\Red{Problems fibred over} 
   {\small $\Purple{W_{\includegraphics[height=5pt]{figures/W3.eps}}}
   \Blue{(X_{\includegraphics[height=5pt]{figures/X1.eps}})^2}
   \ForestGreen{(Y_{\includegraphics[height=5pt]{figures/Y1.eps}})^4}$}\rule{0pt}{25pt}\\
  &
  \begin{tabular}{||l|l|l||}\hline
   Variety& Problem & $\!$\#$\!$\\\hline
  {\small\Blue{\underline{Fl(1,2,3,4;5)}}}\raisebox{-7pt}{\rule{0pt}{19pt}}&
       {\small $\Purple{W_{\includegraphics[height=5pt]{figures/W3.eps}}}
   \Blue{(X_{\includegraphics[height=5pt]{figures/X1.eps}})^2}
   \ForestGreen{(Y_{\includegraphics[height=5pt]{figures/Y1.eps}})^4}
   \Red{Z_{\includegraphics[height=5pt]{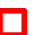}}}$} 
    &7\\\hline
  {\small\Blue{\underline{Fl(1,2,3,4;5)}}}\raisebox{-9pt}{\rule{0pt}{21pt}}&
   {\small $\Purple{W_{\includegraphics[height=5pt]{figures/W3.eps}}}
   \Blue{(X_{\includegraphics[height=5pt]{figures/X1.eps}})^2}
   \ForestGreen{(Y_{\includegraphics[height=5pt]{figures/Y1.eps}})^3}
   \Red{Z_{\includegraphics[height=10pt]{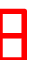}}}$}
    &7\\\hline
  {\small\Blue{\underline{Fl(1,2,3,4;5)}}}\raisebox{-7pt}{\rule{0pt}{19pt}}&
   {\small $A_{4125}
   \Blue{(X_{\includegraphics[height=5pt]{figures/X1.eps}})^2}
   \ForestGreen{(Y_{\includegraphics[height=5pt]{figures/Y1.eps}})^4}$}
    &7\\\hline
  \end{tabular}\vspace{28pt}

 \end{tabular}

\Blue{\rule{430pt}{0.5pt}}

\begin{tabular}{ccc}
    & &{\bf \Red{Point Selection}}\raisebox{-5pt}{\rule{0pt}{10pt}}\\
 \begin{tabular}{||c|c|c|c||}\hline
   \multicolumn{4}{||c||}{\bf Key\rule{0pt}{13pt}}\\\hline
  {\bf Condition}&{\bf Name}&{\bf Symbol}&{\bf Codimension}\\\hline
   \Purple{412}&$\Purple{W_{\includegraphics[height=5pt]{figures/W3.eps}}}$&
   {\bf\Purple{a}}&\Purple{3}\\\hline
   \Blue{132}&$\Blue{X_{\includegraphics[height=5pt]{figures/X1.eps}}}$&
   {\bf\Blue{b}}&\Blue{1}\\\hline
   \ForestGreen{124}&$\ForestGreen{Y_{\includegraphics[height=5pt]{figures/Y1.eps}}}$&
   {\bf \ForestGreen{c}}&\ForestGreen{1}\\\hline
 \end{tabular}
 && \includegraphics[height=40pt]{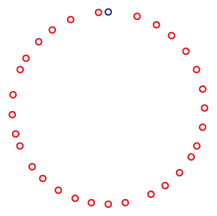}\rule{0pt}{45pt}
\end{tabular}\vspace{5pt}

\Blue{\rule{430pt}{0.5pt}}

Total time of computation: 27,491.26 GHz-seconds or 7.64 GHz-hours on
Noether

\Blue{\rule{430pt}{0.5pt}}

225 000 Polynomial systems solved

\Blue{\rule{430pt}{0.5pt}}

The coefficients of a typical eliminant had 29 digits.

The typical eliminant had size 271 bytes.

\Blue{\rule{430pt}{0.5pt}}

This table automatically generated from the data in 
\Blue{\underline{This File}} using \Blue{\underline{This Maple Script}} 

\Blue{\rule{430pt}{0.5pt}}

\Violet{{\sl Created: Fri Jul 15 15:42:38 CDT 2005}}\raisebox{-5pt}{\rule{0pt}{10pt}}

\end{flushleft}
\end{minipage}
}\end{center}
\caption{Web page for the problem
  $(\Purple{412},(\Blue{132})^2,(\ForestGreen{124})^4)$ on $\Fl(1,2,3;5)$.} 
\label{Fig:html_page}
\end{figure}

After the computations were completed for a given Schubert problem, the data
file was used to generate a web page which displayed information from the
experimentation on that Schubert problem.
Figure~\ref{Fig:html_page} illustrates a typical such page.
This web page has a key in the form of a table with one row for each
Schubert condition.
Each row shows the condition as a permutation, then in a shorthand that is
well-suited to Grassmannian conditions---the letter indicates which member of
the flag it is imposed upon and the partition index indicates the corresponding
Schubert condition on the Grassmannian.
Next is the symbol for that condition used when listing the necklaces, and
finally its codimension.
The figure under \Red{Point Selection} shows the positions of the points in $S$
on $\RP^1$, represented as a circle where the point at the top is $\infty$.
This web page also records the total computation time, the machine used 
(Noether is a computer owned by Sarah Witherspoon), the total number of polynomial
systems solved, and the size of a typical eliminant.

This web page is linked to pages for the problems we
computed which were fibred over it or over which it was
fibred (called \Red{Projections}).
It is also linked to the data file and to the Maple script used
to generate the web page.
At its top is a link \raisebox{-2pt}{\includegraphics[height=12pt]{figures/up.eps}} to 
the web page for the flag variety $\Fl(1,2,3;5)$.
That page lists all 163 Schubert problems we studied on 
$\Fl(1,2,3;5)$, is linked to the other 8 flag varieties in
5-space that we investigated, and to a page with
information about the 29 different flag varieties in our
investigation.

This archive of our data is part of a web page containing 
additional information about this project, which is found at
{\tt www.math.tamu.edu/\~{}sottile/pages/Flags/}, the page displayed 
has further extension {\tt Data/F1235/W3Xe2Ye4.7.html}.
Subsequent addresses will give only the extension from
{\tt .../Flags/Data/}.

\subsection{Design of experiments}
While we investigated many Schubert problems on many small flag
manifolds, by no means did we study all Schubert problems on these
flag manifolds.
We did investigate all Schubert problems on the manifolds of
flags in $\C^4$, and all with degree at least 3 on $\Fl(1,2,3;5)$,
$\Fl(1,2,4;5)$, $\Fl(1,2;5)$, $\Fl(1,3;5)$, $\Fl(2,3;5)$, $\Fl(2,4;5)$,
$\Fl(3,4;5)$, and $\Gr(3,6)$.
Only a small fraction of feasible Schubert problems were investigated
on the other 18 flag manifolds.

There were limitations of
resources which made choices necessary.
For example, the complexity of Gr\"obner basis computation limited
us to Schubert problems of low degree (typically fewer
than 20 solutions).
For the computations on Grassmannians, a more advantageous choice of
local coordinates was possible, which 
allowed significantly larger problems---we studied one
problem on $\Gr(3,7)$ with 91 
solutions\footnote{{\tt F37/We7W2W21.91.html}}.

Many Schubert problems had literally thousands of necklaces, such as
the problem of Table~\ref{table:11352} with 11,352 neckalces.
A systematic study of all necklaces for such a problem would be 
infeasible and the data would be incomprehensible.
We did consider all 1272 necklaces for one such 
problem\footnote{{\tt F12456/Ve2We2W32Ye3Ze2.4.html}}.
Limiting our investigation to problems of small degree and with few necklaces
would still have been infeasible, as there are many thousands of such smaller
Schubert problems on some of these flag manifolds.

On the flag manifolds for which it was impossible to investigate all
Schubert problems, we studied most feasible Grassmannian
Schubert problems, as well as many related to these 
Grassmannian problems through projection, lifting, fibration, and the
notion of child problems as discussed in Sections~\ref{S:T_simple}
and~\ref{S:projections}. 
We looked at some with potentially interesting geometry such as the
problem of Section~\ref{S:Galois}.
We also selected many problems completely at random, intending to sample the
range of possibilities.

Table~\ref{T:X} lists the Schubert problems discussed here,
their associated web pages, and the resources expended on each.

\begin{table}[htb]
\begin{tabular}{|c|l|r|}\hline
 Location& \multicolumn{1}{c|}{Web Page} & \multicolumn{1}{c|}{CPU} \\\hline
Table~\ref{table:12-flag}& {\tt F235/Xe4Ye4.12.html}
       &213.38 GHz-days\\\hline
Table~\ref{table:new12-flag}&{\tt F12345/We2Xe3Ye3Ze2.12.html}
       & 47.61 GHz-days\\\hline
Table~\ref{table:11352}&{\tt F123456/Ve2We2XX321Ye2Ze2.8.html}
       & 5.25 GHz-days \\\hline
Table~\ref{table:lower}&{\tt F347/WW31e2Xe2X211.10.html}
       &  63.43 GHz-days \\\hline
Table~\ref{T:upper}&{\tt F12345/A1325e2A2143e3.7.html}
       & 1.94 GHz-days \\\hline
Table~\ref{T:gaps}&{\tt F1356/A21436e2A31526Xe2.8.html}
       & 12.84 GHz-days \\\hline
Table~\ref{T:seven}&{\tt F2346/A1432A1254We2X21Ye2.12.html}
       & 61.86 GHz-days\\\hline
Table~\ref{T:fibration}&{\tt  F246/We4Y11e4.6.html}
       &  13.57 GHz-days \\\hline
Figure~\ref{Fig:table_discr}&{\tt F12456/A13254A43125A12564VZ.2.html}
       & 1.31 GHz-days\\\hline

\end{tabular}\vspace{5pt}
\caption{CPU time used for computations shown here\label{T:X}}
\end{table}

\section{Conclusion and Future Work}
We presented a geometrically vivid example of the failure of the Shapiro
Conjecture for Schubert intersections given by osculating flags on flag
manifolds, and presented a refinement of the conjecture for flag varieties. 
Significant evidence, both theoretical and experimental, was presented in
support of this refinement.
Several new phenomena discovered in this experimentation were
presented.

The proof of the conjecture for certain two-step flag manifolds by Eremenko 
{\it  et al.} leads to an extension concerning secant flags.
The further investigation of this secant flag conjecture is a worthwhile future
project. 


\providecommand{\bysame}{\leavevmode\hbox to3em{\hrulefill}\thinspace}
\providecommand{\MR}{\relax\ifhmode\unskip\space\fi MR }
\providecommand{\MRhref}[2]{%
  \href{http://www.ams.org/mathscinet-getitem?mr=#1}{#2}
}
\providecommand{\href}[2]{#2}


\begin{thebibliography}{10}

\bibitem{EH83}
D.~Eisenbud and J.~Harris, \emph{Divisors on general curves and cuspidal
  rational curves}, Invent. Math. \textbf{74} (1983), 371--418.

\bibitem{EH87}
\bysame, \emph{When ramification points meet}, Invent. Math. \textbf{87}
  (1987), 485--493.

\bibitem{EG01b}
A.~Eremenko and A.~Gabrielov, \emph{Degrees of real {W}ronski maps}, Discrete
  Comput. Geom. \textbf{28} (2002), no.~3, 331--347.

\bibitem{EG02a}
\bysame, \emph{Rational functions with real critical points and the {B}. and
  {M}. {S}hapiro conjecture in real enumerative geometry}, Ann. of Math. (2)
  \textbf{155} (2002), no.~1, 105--129.

\bibitem{EGSV}
A.~Eremenko, A.~Gabrielov, M.~Shapiro, and A.~Vainshtein, \emph{{Rational
  functions and real Schubert calculus}}, Proc.~AMS, to appear.

\bibitem{FGLM}
J.-C. Faug\`ere, P.~Gianni, D.~Lazard, and T.~Mora, \emph{Efficient computation
  of zero-dimensional {G}roebner bases by change of ordering}, J. Symb. Comp.
  \textbf{16} (1993), 329--344.

\bibitem{Fu97}
W.~Fulton, \emph{Young tableau}, Cambridge University Press, 1997.

\bibitem{SINGULAR}
G.-M. Greuel, G.~Pfister, and H.~Sch\"onemann, \emph{{\sc Singular} 2.0}, {A
  Computer Algebra System for Polynomial Computations}, Centre for Computer
  Algebra, University of Kaiserslautern, 2001, {\tt
  http://www.singular.uni-kl.de}.

\bibitem{Ha79}
J.~Harris, \emph{Galois groups of enumerative problems}, Duke Math.~J.
  \textbf{46} (1979), 685--724.

\bibitem{IKS}
I.~V. Itenberg, V.~M. Kharlamov, and E.~I. Shustin, \emph{Logarithmic
  equivalence of the {W}elschinger and the {G}romov-{W}itten invariants},
  Uspekhi Mat. Nauk \textbf{59} (2004), no.~6(360), 85--110.

\bibitem{KhS03}
V.~Kharlamov and F.~Sottile, \emph{Maximally inflected real rational curves},
  Moscow Math. J., {\bf 3}, 2003.

\bibitem{Kl74}
S.~Kleiman, \emph{The transversality of a general translate}, Compositio Math.
  \textbf{28} (1974), 287--297.

\bibitem{Mi}
Grigory Mikhalkin, \emph{Enumerative tropical algebraic geometry in
  {$\mathbb{R}^2$}}, J. Amer. Math. Soc. \textbf{18} (2005), no.~2, 313--377
  (electronic).

\bibitem{MV04}
E.~Mukhin and A.~Varchenko, \emph{Critical points of master functions and flag
  varieties}, Commun. Contemp. Math. \textbf{6} (2004), no.~1, 111--163.

\bibitem{PuSo}
Kevin Purbhoo and Frank Sottile, \emph{The recursive nature of the cominuscule
  {S}chubert calculus}, 2005, Manuscript in progress.

\bibitem{RS98}
J.~Rosenthal and F.~Sottile, \emph{Some remarks on real and complex output
  feedback}, Systems \& Control Lett. \textbf{33} (1998), no.~2, 73--80, For a
  description of the computational aspects, see {\tt
  www.math.tamu.edu/\~{}sottile/pages/control/}.

\bibitem{Sch00}
Claus Scheiderer, \emph{Sums of squares of regular functions on real algebraic
  varieties}, Trans. Amer. Math. Soc. \textbf{352} (2000), no.~3, 1039--1069.

\bibitem{SS02}
V.~Sedykh and B~Shapiro, \emph{Two conjectures on convex curves},
Intern.~J.~Math., to appear. 

\bibitem{SS}
E.~Soprunova and F.~Sottile, \emph{Lower bounds for real solutions to sparse
  polynomial systems}, Adv.~Math., to appear.

\bibitem{So96}
F.~Sottile, \emph{Pieri's formula for flag manifolds and schubert polynomials},
  Ann.~Inst.~Fourier \textbf{46} (1996), 1--22.

\bibitem{So97c}
\bysame, \emph{Enumerative geometry for real varieties}, Algebraic Geometry,
  Santa Cruz 1995 (J.~Koll{\'a}r, R.~Lazarsfeld, and D.~Morrison, eds.),
  Proc.~Sympos.~Pure Math., vol. 62, Part 1, Amer.~Math.~Soc., 1997,
  pp.~435--447.

\bibitem{So97a}
\bysame, \emph{Enumerative geometry for the real {G}rassmannian of lines in
  projective space}, Duke Math. J. \textbf{87} (1997), no.~1, 59--85.

\bibitem{So99a}
\bysame, \emph{The special {S}chubert calculus is real}, ERA of the AMS
  \textbf{5} (1999), 35--39.

\bibitem{So00a}
\bysame, \emph{Real {S}chubert calculus: Polynomial systems and a conjecture of
  {S}hapiro and {S}hapiro}, Exper.~Math. \textbf{9} (2000), 161--182.

\bibitem{So00b}
\bysame, \emph{Some real and unreal enumerative geometry for flag manifolds},
  Mich. Math. J. \textbf{48} (2000), 573--592, Special Issue in Honor of
  Wm.~Fulton.

\bibitem{So_M2}
\bysame, \emph{From enumerative geometry to solving systems of
  polynomials equations}, Computations in algebraic geometry with Macaulay 2,
  Algorithms Comput. Math., vol.~8, Springer, Berlin, 2002, pp.~101--129.

\bibitem{So03b}
\bysame, \emph{Enumerative real algebraic geometry}, Algorithmic and
  quantitative real algebraic geometry (Piscataway, NJ, 2001), DIMACS Ser.
  Discrete Math. Theoret. Comput. Sci., vol.~60, Amer. Math. Soc., Providence,
  RI, 2003, pp.~139--179.

\bibitem{Va04}
R.~Vakil, \emph{{Schubert induction}}, Annals of Math., to appear.

\bibitem{Ver00}
J.~Verschelde, \emph{Numerical evidence of a conjecture in real algebraic
  geometry}, Exper.~Math. \textbf{9} (2000), 183--196.

\bibitem{W}
Jean-Yves Welschinger, \emph{Invariants of real rational symplectic 4-manifolds
  and lower bounds in real enumerative geometry}, C. R. Math. Acad. Sci. Paris
  \textbf{336} (2003), no.~4, 341--344.

\end{thebibliography}
\end{document}